\documentclass[11pt]{article}
\usepackage{amssymb,a4wide,latexsym,makeidx,epsfig}
\usepackage{amsthm}
\usepackage{dsfont}
\usepackage{amsmath}
\usepackage{lipsum}
\usepackage{enumerate}
\usepackage{mathrsfs}
\usepackage{xcolor}
\usepackage{tikz}
\usepackage{setspace}
\usepackage{geometry}
\usepackage{appendix}
\usepackage{multirow}
\usepackage{subcaption}
\usepackage[colorlinks=true, linkcolor=black, anchorcolor=black, citecolor=black, urlcolor=black, CJKbookmarks=true]{hyperref}
\usepackage{microtype}
\usepackage{colonequals}
\allowdisplaybreaks

\AtBeginDocument{
  \setlength{\abovedisplayskip}{6pt} %公式上方间距（长公式，前一行文字较长时）
  \setlength{\belowdisplayskip}{6pt} %公式下方间距
  \setlength{\abovedisplayshortskip}{0pt} %公式上方间距（短公式，前一行文字较短时）
  \setlength{\belowdisplayshortskip}{3pt} %公式下方间距（短公式）
  \setlength{\jot}{4pt}
}

\newtheorem{theorem}{Theorem}[section]

\newtheorem{lemma}[theorem]{Lemma}

\newtheorem{observation}[theorem]{Observation}
\newtheorem{proposition}[theorem]{Proposition}
\newtheorem{corollary}[theorem]{Corollary}
\newtheorem{problem}[theorem]{Problem}

\newtheorem{claim}{Claim}[section]

\numberwithin{figure}{section}
\numberwithin{equation}{section}

\begin{document}

\title{Edge-disjoint Hamilton cycles in the McKay--Miller--\v{S}ir\'{a}\v{n} graphs $H_q$ for every prime power $q\geq 3$}

\author{
%Xihe Li\footnote{School of Mathematics and Statistics, Shaanxi Normal University, Xi'an, Shaanxi 710119, China. E-mail address: xiheli@snnu.edu.cn.}
Xihe Li \vspace{0.2cm}\\
\small{School of Mathematics and Statistics, Shaanxi Normal University, Xi'an, Shaanxi 710119, China.} \\
\small{xiheli@snnu.edu.cn}
}
\date{}
\maketitle
\newcommand\blfootnote[1]{%
\begingroup
\renewcommand\thefootnote{}\footnote{#1}%
\addtocounter{footnote}{-1}%
\endgroup
}
\blfootnote{This work was supported by National Natural Science Foundation of China (Grant No. 12501492), Shaanxi Fundamental Science Research Project for Mathematics and Physics (Grant No. 25JSQ043), and Fundamental Research Funds for the Central Universities (Grant No. GK202506024).}
\begin{center}
\begin{minipage}{145mm}
\vskip -0.3cm
\begin{center}
{\small {\bf Abstract}}
\end{center}
{\small Motivated by Hamilton cycle problems and applications to collective communication, we study edge-disjoint Hamilton cycles in the McKay--Miller--\v{S}ir\'{a}\v{n} (MMS) graphs.
The MMS graph $H_q$ (indexed by a prime power $q$) is a highly symmetric graph arising from finite fields.
For every prime power $q\geq 3$, we prove that the MMS graph $H_q$ contains $\big\lceil\frac{q-1}{2}\big\rceil$ edge-disjoint Hamilton cycles.
This resolves a problem of Li, Wang and Zheng, who established the result for prime $q$ and for $q=4$.
Our construction uses graph lifts over finite fields and joins the components of lifted Hamilton cycles by adding edges within vertex fibres.
The resulting cycles can be constructed using $O(|E(H_q)|)$ finite field operations, comparisons, and elementary list operations.
As an application, these cycles provide edge-disjoint spanning rings for parallel ring-all-reduce algorithms on MMS networks.

\hspace{2em}

}
\end{minipage}
\end{center}

\section{Introduction}
\label{sec:introduction}

McKay--Miller--\v{S}ir\'{a}\v{n} graphs form an important family of highly symmetric graphs arising from finite fields.
In the study of the \emph{degree-diameter problem}, McKay, Miller and \v{S}ir\'{a}\v{n}~\cite{McMS} constructed a family of graphs of diameter 2, now known as the \emph{McKay--Miller--\v{S}ir\'{a}\v{n} graphs} (MMS graphs).
Their orders attain asymptotically $\frac{8}{9}$ of the Moore bound for diameter 2.
The MMS graph $H_q$ (indexed by a prime power $q$) can be constructed via a \emph{graph lift operation} on a base graph obtained from the complete bipartite graph $K_{q,q}$ by attaching loops or semiedges.
We postpone the explicit definition to Section~\ref{subsec:H_q}, since several additional concepts concerning graph lifts and finite fields are needed.
For $q\equiv 1 \pmod{4}$, the graph $H_q$ is vertex-transitive.
Moreover, Dinitz, Schapira and Shahaf~\cite{DiSS} proved that MMS graphs have good expansion properties.
In this paper, we study edge-disjoint Hamilton cycles in MMS graphs.
Our motivation comes from Hamilton cycle problems for vertex-transitive graphs and applications to interconnection network design and collective communications.

All graphs considered in this paper are finite.
A \emph{Hamilton cycle} (\emph{path}) in a graph $G$ is a cycle (path) containing all vertices of $G$.
The problem of finding Hamilton cycles or paths is a fundamental research topic in graph theory.
A graph is \emph{vertex-transitive} if its automorphism group acts transitively on its vertices.
In 1969, Lov\'{a}sz~\cite{Lov970} asked whether every connected vertex-transitive graph admits a Hamilton path, which is now known as the Lov\'{a}sz conjecture.
A stronger conjecture of Thomassen (see~\cite{Ber978,KuMa}) states that there are only finitely many connected vertex-transitive graphs without Hamilton cycles, while Babai~\cite{Bab979,Bab995} conjectured that there are infinitely many such graphs.
For Cayley graphs, the Hamilton cycle problem was posed by Rapaport-Strasser~\cite{Rap} in 1959.
It is conjectured that every finite connected Cayley graph with at least three vertices has a Hamilton cycle.
Moreover, Alspach~\cite{Als984} conjectured that every $2k$-regular connected Cayley graph on a finite abelian group can be decomposed into $k$ edge-disjoint Hamilton cycles.
These conjectures have attracted wide attention over the past several decades.
For more information, we refer the interested reader to \cite{BrDe,BCPS26+,ChHM,LiMe,MaZh,MeMN} for results on vertex-transitive graphs, and \cite{AMRS,BeFM,ErLe,ErLP,GlCS,LeMM,Liu03} for results on Cayley graphs.
For Hamilton cycles in lifts of graphs, see \cite{BCCF,ChFr,LuWW}.

In 2014, Besta and Hoefler~\cite{BeHo} introduced the \emph{Slim Fly network} using MMS graphs as the underlying topology.
This topology has attracted significant attention in high-performance computing due to its performance metrics, including latency, bandwidth, cost, resiliency and power consumption; see \cite{BeHo,Blaetal,HaWB,RMFY,WCMRC,YEGQH}.
The \emph{all-reduce} operation is an important collective communication operation in which data from all participating vertices are combined using a reduction operation, and the result is made available at every participating vertex.
The ring-all-reduce algorithm, which makes use of a cycle structure, is a bandwidth-optimal all-reduce algorithm.
If a network contains $t$ edge-disjoint Hamilton cycles, the data at each vertex can be divided into $t$ parts, and a ring-all-reduce algorithm can be run on each cycle simultaneously, with different cycles handling different parts of the data; see \cite{Hun15,WKZGJMZK}.
This motivates the construction of many edge-disjoint Hamilton cycles in MMS graphs.

It is well-known that if $q\geq 3$ is odd, then $K_{q,q}$ can be decomposed into $\frac{q-1}{2}$ edge-disjoint Hamilton cycles and a perfect matching;
if $q\geq 2$ is even, then $K_{q,q}$ can be decomposed into $\frac{q}{2}$ edge-disjoint Hamilton cycles.
We will present the decomposition in Section~\ref{subsec:cycles}.
In~\cite{LiWZ}, Li, Wang and Zheng proposed the following problem.

\begin{problem}[\cite{LiWZ}]\label{prob:q}
For any prime power $q\geq 3$, does the MMS graph $H_{q}$ contain $\big\lceil\frac{q-1}{2}\big\rceil$ edge-disjoint Hamilton cycles?
\end{problem}

They answered Problem~\ref{prob:q} affirmatively for every prime $q\geq 3$ and for $q=4$.
In this paper, we solve this problem completely.

\begin{theorem}\label{thm:q}
For any prime power $q\geq 3$, the MMS graph $H_{q}$ contains $\big\lceil\frac{q-1}{2}\big\rceil$ edge-disjoint Hamilton cycles.
\end{theorem}

For completeness, we will provide a full proof of Theorem~\ref{thm:q}, including the cases of prime $q$ and $q=4$ established in~\cite{LiWZ}.
Our proof is essentially constructive, and the $\big\lceil\frac{q-1}{2}\big\rceil$ edge-disjoint Hamilton cycles given in Theorem~\ref{thm:q} can be constructed using $O(|E(H_q)|)$ finite field operations, comparisons, and elementary list operations.
We will prove this in Appendix~\ref{sec:app-2}.
As an application, for networks with MMS graphs as underlying graphs, the resulting cycles provide edge-disjoint spanning rings that can carry separate portions of a parallel ring-all-reduce operation.

The remainder of this paper is organized as follows.
In the next section, we present the definition of the MMS graphs $H_q$, describe the lifts in $H_q$ of Hamilton cycles in $K_{q,q}$, and provide a proof sketch of our result.
In Section~\ref{sec:join}, we show how to join the components of a lifted Hamilton cycle into a Hamilton cycle.
In Section~\ref{sec:proof}, we finish the proof of Theorem~\ref{thm:q}.
Appendix~\ref{sec:app-2} establishes the operation count for the construction, and Appendix~\ref{sec:app-1} illustrates the graphs for $q=3,4$.

\section{Preliminaries}
\label{sec:pre}

In this section, we first present several elementary facts on finite fields and provide the definition of the MMS graphs in Section~\ref{subsec:H_q}.
There are three equivalent definitions for the MMS graphs: lift of a modified $K_{q,q}$~\cite{McMS}, lift of dipoles~\cite{Sia}, and a geometric construction~\cite{Haf}.
In order to prove our main result, we choose to present the definition as stated in~\cite{McMS} (see also~\cite{LiWZ}).
In Section~\ref{subsec:cycles}, we construct $\big\lceil\frac{q-1}{2}\big\rceil$ edge-disjoint Hamilton cycles in $K_{q,q}$, and describe their lifts in $H_q$.
In Section~\ref{subsec:sketch}, we provide a proof sketch of Theorem~\ref{thm:q}.

\subsection{Finite fields and the MMS graph $H_q$}
\label{subsec:H_q}

We first recall some elementary facts on finite fields.
Let $q=p^n$ be a prime power with $q\geq 3$, where $p$ is prime and $n$ is a positive integer.
Then $q\equiv -1, 0, 1 \pmod{4}$.
We write $q=4\ell+\delta$, where $\ell$ is a positive integer and $\delta \in \{-1, 0, 1\}$.
Let $\mathbb{F}_q$ be the finite field of order $q$ with addition operation $+$ and multiplication operation $\cdot$, and let $\mathbb{F}_q^+$ be its additive group and $\mathbb{F}_q^{\times}$ be its multiplicative group.
The field $\mathbb{F}_q$ contains $q$ elements and has characteristic $p$.
In particular, for any element $g\in \mathbb{F}_q$, we have $pg=0$.
The elements $0,1,\ldots,p-1$ form a subfield $\mathbb{F}_p$ of $\mathbb{F}_q$.
With addition and scalar multiplication inherited from the field operations, $\mathbb{F}_q$ is an $n$-dimensional vector space over $\mathbb{F}_p$.
In particular, if $a\neq 0$, then $a\mathbb{F}_p=\{0,a,2a,\ldots,(p-1)a\}$ has exactly $p$ elements.
Let $\xi$ be a primitive element of $\mathbb{F}_q$, so $\xi^{q-1}=1$ and $\mathbb{F}_q=\{\xi, \xi^2, \ldots, \xi^{q-2}, \xi^{q-1}, 0\}$.
If $q$ is odd, then $\xi^{(q-1)/2}=-1$; if $q$ is even, then $1=-1$.

We next introduce the \emph{graph lift operation}, which has many applications in graph theory \cite{BiLi,BoCo,MaSS,Pud} and other areas \cite{KLVW,PaKa,SmMi}.
Let $G=(V(G), E(G))$ be an undirected graph, possibly with loops, parallel edges and semiedges, where a semiedge is a dangling edge with just one end incident to a vertex.
Each edge (including loops) other than a semiedge gives rise to two opposite arcs.
A semiedge has only one possible direction, i.e., outward from the incident vertex.
The \emph{arc set} $D(G)$ of $G$ is the set of all arcs of $G$; the reverse of an arc $e$ is denoted by $e^{-1}$.
If $e$ is a semiedge, we may still use the symbol $e^{-1}$, but we understand that $e=e^{-1}$ in this case.
Let $\Gamma$ be a finite group with binary operation $\circ$ and identity element $e_{\Gamma}$.
A mapping $\alpha: D(G)\to \Gamma$ is a \emph{voltage assignment} if
\begin{equation}\label{eq:voltage}
\alpha(e^{-1})=(\alpha(e))^{-1}
\end{equation}
for each arc $e\in D(G)$, where $(\alpha(e))^{-1}$ is the inverse element of $\alpha(e)$ in $\Gamma$.

Given $G$, $\Gamma$ and a voltage assignment $\alpha$, a \emph{lift} of $G$ is a new graph, denoted by $G^{\alpha}$, with vertex set $V(G)\times \Gamma$ and edge set defined as follows.
For an arc $e\in D(G)$ from vertex $u$ to vertex $v$ that is not a semiedge, and for each $g\in\Gamma$, the lift has an arc $(e,g)$ from $(u,g)$ to $(v,g\circ\alpha(e))$.
By definition, the arc set of $G^{\alpha}$ lifted from $D(G)$ is $D(G)\times \Gamma$.
Moreover, the arcs $(e,g)$ and $(e^{-1}, g\circ \alpha(e))$ are mutually reverse in the lift, and we view them as an undirected edge $e_g$ of $G^{\alpha}$, so the lift is an undirected graph.
For a semiedge $e$ at some vertex $u\in V(G)$, we have $\alpha(e)\circ \alpha(e)=e_{\Gamma}$ by \eqref{eq:voltage}.
If $\alpha(e)=e_{\Gamma}$, then $e$ lifts to $|\Gamma|$ semiedges;
if $\alpha(e)$ is a nontrivial element of order 2 in $\Gamma$, then $g\circ \alpha(e)\circ \alpha(e)=g$, and $e$ lifts to a set $\big\{(u,g)(u,g\circ \alpha(e))\colon\, g\in \Gamma\big\}$ of $\frac{|\Gamma|}{2}$ ordinary edges, rather than semiedges.

The graph $G$ is called the \emph{base graph} of $G^{\alpha}$.
For each vertex $v\in V(G)$, the set $\{(v, g)\colon\, g\in \Gamma\}$ of vertices in $G^{\alpha}$ is called the \emph{fibre} above $v$;
for each arc $e\in D(G)$, the set $\{(e, g)\colon\, g\in \Gamma\}$ of arcs in $G^{\alpha}$ is called the \emph{fibre} above $e$.
We use $\phi(v)$ and $\phi(e)$ to denote the fibre above a vertex $v$ and an arc $e$, respectively.
For an undirected edge $e$, we also write $\phi(e)$ for the set of undirected edges lying above $e$.
When a voltage is used, an orientation of $e$ is specified.

Now we have all the ingredients to present an explicit definition of the MMS graphs $H_q$.
The base graph $G_q$ is obtained from the complete bipartite graph $K_{q,q}$ by attaching $\ell$ loops at each vertex if $q$ is odd, and attaching $2\ell=\frac{q}{2}$ semiedges at each vertex if $q$ is even.
Let $U\colonequals \{u_x\colon\, x\in \mathbb{F}_q\}$ and $V\colonequals \{v_y\colon\, y\in \mathbb{F}_q\}$ be the partite sets of $K_{q,q}$.
Recall that $q=4\ell+\delta$ with $\delta \in \{-1, 0, 1\}$.
For odd $q$, let $Y=\{\xi^{2j}:0\leq j\leq \ell-1\}$.
If $\delta=1$, let $Y'=\{\xi^{2j+1}:0\leq j\leq \ell-1\}$; if $\delta=-1$, let $Y'=\{1\}\cup\{\xi^{2j-1}:1\leq j\leq \ell-1\}$.
For even $q$, let $X=\{\xi^{2j}:0\leq j\leq 2\ell-1\}$ and $X'=\{\xi^{2j+1}:0\leq j\leq 2\ell-1\}$.
The voltage assignment $\alpha: D(G_q) \to \mathbb{F}_q^+$ is defined as follows:
\begin{itemize}
\item for all $q$, each arc from a vertex $u_x$ to a vertex $v_y$ will receive $xy$ as its voltage, and each arc from a vertex $v_y$ to a vertex $u_x$ will receive $-xy$ as its voltage, where $xy=x \cdot y$ is the multiplication in $\mathbb{F}_q$;
\item for odd $q$, $\alpha$ assigns to the $\ell$ loops\footnote{For every loop, choose one of its two orientations. The specified voltage is assigned to this orientation, and its negative is assigned to the reverse orientation.} at each vertex $u_x$ the voltages in $Y$ bijectively, and $\alpha$ assigns to the $\ell$ loops at each vertex $v_y$ the voltages in $Y'$ bijectively;
\item for even $q$, $\alpha$ assigns to the $2\ell$ semiedges at each vertex $u_x$ the voltages in $X$ bijectively, and $\alpha$ assigns to the $2\ell$ semiedges at each vertex $v_y$ the voltages in $X'$ bijectively.
\end{itemize}
The \emph{MMS graph} $H_q\colonequals (G_q)^{\alpha}$ is the lift of $G_q$ with respect to the voltage assignment $\alpha$ and group $\mathbb{F}_q^+$.

For clarity, we now describe the lift explicitly (see Figures~\ref{fig:q=3} and \ref{fig:q=4} in Appendix~\ref{sec:app-1} for an illustration).
Every vertex $u_x\in U$ (resp., $v_y\in V$) lifts to a set $\{(u_{x},g)\colon\, g\in \mathbb{F}_q\}$ (resp., $\{(v_y, g)\colon\, g\in \mathbb{F}_q\}$) of $q$ vertices of $H_q$.
Every edge $u_xv_y$ of $K_{q,q}$ lifts to a matching $\{(u_x, g)(v_{y},g+xy)\colon\, g\in \mathbb{F}_q\}$ of $q$ edges of $H_q$.
We call these edges the \emph{cross edges}.
If $q$ is odd (i.e., $\delta\in \{-1,1\}$), then a loop of voltage $s\in Y$ at $u_x\in U$ lifts to a set $\{(u_x, g)(u_{x},g+s)\colon\, g\in \mathbb{F}_q\}$ of edges in $H_q$, and a loop of voltage $s\in Y'$ at $v_y\in V$ lifts to a set $\{(v_y, g)(v_{y},g+s)\colon\, g\in \mathbb{F}_q\}$ of edges in $H_q$.
If $q$ is even (i.e., $\delta=0$), then $q=2^n$.
Thus every element $s\in \mathbb{F}_q$ satisfies $s+s=0$, so $g+s+s=g$ for any $g\in \mathbb{F}_q$.
A semiedge of voltage $s\in X$ at $u_x\in U$ lifts to a matching $\{(u_x, g)(u_{x},g+s)\colon\, g\in \mathbb{F}_q\}$ of size $\frac{q}{2}$ in $H_q$, and a semiedge of voltage $s\in X'$ at $v_y\in V$ lifts to a matching $\{(v_y, g)(v_{y},g+s)\colon\, g\in \mathbb{F}_q\}$ of size $\frac{q}{2}$ in $H_q$.

Let $S_U\colonequals Y\cup (-Y)$ if $q$ is odd, and $S_U\colonequals X$ if $q$ is even.
Let $S_V\colonequals Y'\cup (-Y')$ if $q$ is odd, and $S_V\colonequals X'$ if $q$ is even.
Then $S_U=-S_U$, $S_V=-S_V$ and $|S_U|=|S_V|=\frac{q-\delta}{2}$.
Moreover, the set of edges inside the fibre above a vertex $u_x\in U$ is exactly $\{(u_x, g)(u_{x},g+s)\colon\, g\in \mathbb{F}_q, s\in S_U\}$, and the set of edges inside the fibre above a vertex $v_y\in V$ is exactly $\{(v_y, g)(v_{y},g+s)\colon\, g\in \mathbb{F}_q, s\in S_V\}$.
We call $S_U$ (resp., $S_V$) the \emph{difference set} of the fibre above a vertex of $U$ (resp., $V$).
It follows that $H_q$ has $2q^2$ vertices; each vertex is incident with $q$ cross edges and $|S_U|$ or $|S_V|$ edges inside a fibre.
In particular, $H_q$ is a $\frac{3q-\delta}{2}$-regular graph.

\begin{lemma}\label{le:difference}
We have $S_U\cup S_V=\mathbb{F}_q^{\times}$, and each of $S_U$ and $S_V$ spans $\mathbb{F}_q$ as a vector space over $\mathbb{F}_p$.
Moreover, $S_U\cap S_V=\varnothing$ when $\delta=1$, $S_U\cap S_V=\{1,-1\}$ when $\delta=-1$, and $S_U\cap S_V=\{1\}$ when $\delta=0$.
\end{lemma}

\begin{proof}
If $\delta=1$, then $\xi^{4\ell}=\xi^{q-1}=1$ and $\xi^{2\ell}=\xi^{(q-1)/2}=-1$.
Thus
\begin{align*}
S_U = &~Y\cup (-Y)=\big\{1, \xi^2, \xi^4, \ldots, \xi^{2\ell-2}\big\}\cup \big\{-1, -\xi^2, -\xi^4, \ldots, -\xi^{2\ell-2}\big\} \\
= &~\big\{1, \xi^2, \xi^4, \ldots, \xi^{2\ell-2}\big\}\cup \big\{\xi^{2\ell}, \xi^{2\ell+2}, \xi^{2\ell+4}, \ldots, \xi^{4\ell-2}\big\}=\big\{1, \xi^2, \xi^4, \ldots, \xi^{4\ell-2}\big\}, \\
S_V = &~Y'\cup (-Y')=\big\{\xi, \xi^3, \ldots, \xi^{2\ell-1}\big\}\cup \big\{-\xi, -\xi^3, \ldots, -\xi^{2\ell-1}\big\} \\
= &~\big\{\xi, \xi^3, \ldots, \xi^{2\ell-1}\big\}\cup \big\{\xi^{2\ell+1}, \xi^{2\ell+3}, \ldots, \xi^{4\ell-1}\big\}=\big\{\xi, \xi^3, \ldots, \xi^{4\ell-1}\big\}.
\end{align*}
Hence, $S_U\cup S_V=\mathbb{F}_q^{\times}$ and $S_U\cap S_V=\varnothing$.

If $\delta=-1$, then $\xi^{4\ell-2}=\xi^{q-1}=1$ and $\xi^{2\ell-1}=\xi^{(q-1)/2}=-1$.
Thus
\begin{align*}
S_U = &~Y\cup (-Y)=\big\{1, \xi^2, \xi^4, \ldots, \xi^{2\ell-2}\big\}\cup \big\{-1, -\xi^2, -\xi^4, \ldots, -\xi^{2\ell-2}\big\} \\
= &~\big\{1, \xi^2, \xi^4, \ldots, \xi^{2\ell-2}\big\}\cup \big\{-1, \xi^{2\ell+1}, \xi^{2\ell+3}, \ldots, \xi^{4\ell-3}\big\}, \\
S_V = &~Y'\cup (-Y')=\big\{1, \xi, \xi^3, \ldots, \xi^{2\ell-3}\big\}\cup \big\{-1, -\xi, -\xi^3, \ldots, -\xi^{2\ell-3}\big\} \\
= &~\big\{1, \xi, \xi^3, \ldots, \xi^{2\ell-3}\big\}\cup \big\{-1, \xi^{2\ell}, \xi^{2\ell+2}, \ldots, \xi^{4\ell-4}\big\}.
\end{align*}
Hence, $S_U\cup S_V=\mathbb{F}_q^{\times}$ and $S_U\cap S_V=\{1,-1\}$.

If $\delta=0$, then $\mathbb{F}_q$ has characteristic $2$, so $-1=1=\xi^{q-1}=\xi^{4\ell-1}$.
Thus $S_U= X=\{1, \xi^2, \xi^4, \ldots, \xi^{4\ell-2}\}$ and $S_V= X'=\{\xi, \xi^3, \ldots, \xi^{4\ell-1}\}$ satisfy $S_U\cup S_V=\mathbb{F}_q^{\times}$ and $S_U\cap S_V=\{1\}$.

We next show that each of $S_U$ and $S_V$ spans $\mathbb{F}_q$ as a vector space over $\mathbb{F}_p$.
Recall that $\mathbb{F}_q$ is an $n$-dimensional vector space over $\mathbb{F}_p$.
A proper subspace of $\mathbb{F}_q$ has at most $p^{n-1}=\frac{q}{p}$ elements, and thus contains at most $\frac{q}{p}-1$ nonzero elements.
Note that neither $S_U$ nor $S_V$ contains $0$.
Moreover, $|S_U|=|S_V|=\frac{q-\delta}{2}$.
If $\delta\in \{-1,1\}$, then $p\geq 3$ is odd and $\frac{q-\delta}{2}\geq \frac{q-1}{2} > \frac{q}{p}-1$;
if $\delta=0$, then $p=2$ and $\frac{q-\delta}{2}= \frac{q}{2} > \frac{q}{2}-1$.
Thus neither $S_U$ nor $S_V$ can lie inside a proper subspace of $\mathbb{F}_q$.
Therefore, each of $S_U$ and $S_V$ spans $\mathbb{F}_q$ as a vector space over $\mathbb{F}_p$.
\end{proof}

\subsection{Hamilton cycles in $K_{q,q}$ and their lifts in $H_q$}
\label{subsec:cycles}

We first construct the $\big\lceil\frac{q-1}{2}\big\rceil$ edge-disjoint Hamilton cycles in the $K_{q,q}$ subgraph of the base graph $G_q$.
Recall that $V(G_q)=U\cup V$, $U= \{u_x\colon\, x\in \mathbb{F}_q\}$, $V= \{v_y\colon\, y\in \mathbb{F}_q\}$, and $\mathbb{F}_q=\{\xi, \xi^2, \ldots, \xi^{q-2}, \xi^{q-1}(=1), 0\}$.
For simplicity, throughout Sections~\ref{subsec:cycles}--\ref{sec:proof} and Appendix~\ref{sec:app-2}, we write $u_i\colonequals u_{\xi^{i}}$ and $v_i\colonequals v_{\xi^{i}}$ for each $i\in \{1, 2, \ldots, q-1\}$, and write $u_q\colonequals u_0$ and $v_q\colonequals v_0$.
In particular, under this convention, $u_1$ is $u_{\xi}$, not $u_{\xi^{q-1}}$, and $v_1$ is $v_{\xi}$, not $v_{\xi^{q-1}}$.
Now $U=\{u_1, u_2, \ldots, u_q\}$ and $V=\{v_1, v_2, \ldots, v_q\}$.
Here and below, the indices in $u_i$ and $v_i$ are considered modulo $q$ and represented in $\{1, 2, \ldots, q\}$, whereas exponents of $\xi$ are reduced modulo $q-1$.

For each $s\in \{1, 2, \ldots, q\}$, let $M_s\colonequals \{u_iv_{i+s-1}\colon\, 1\leq i\leq q\}$ be a perfect matching of $K_{q,q}$.
Then $M_1, M_2, \ldots, M_q$ are pairwise edge-disjoint and $\bigcup_{1\leq s\leq q}M_s=E(K_{q,q})$.
For each $k\in \{1, 2, \ldots, \big\lceil\frac{q-1}{2}\big\rceil\}$, let $C^k$ be the subgraph of $K_{q,q}$ consisting of edges in $M_{2k}\cup M_{2k+1}$, where $M_{q+1}=M_1$.
Then each $C^k$ is a Hamilton cycle.
Indeed, starting at $v_{2k}$, we can pass through all vertices as follows:
\begin{equation}\label{eq:Ck}
v_{2k}\to u_1\to v_{2k+1}\to u_2\to v_{2k+2}\to \cdots \to u_q \to v_{2k}.
\end{equation}
At each step, we use alternately an edge of $M_{2k}$ and an edge of $M_{2k+1}$.
Moreover, $C^1, C^2, \ldots, C^{\lceil(q-1)/2\rceil}$ are pairwise edge-disjoint, since they use disjoint pairs of matchings.
In this way, when $q$ is odd, we decompose $E(K_{q,q})$ into $\frac{q-1}{2}$ edge-disjoint Hamilton cycles $C^1, C^2, \ldots, C^{(q-1)/2}$ and a perfect matching $M_1$;
when $q$ is even, we decompose $E(K_{q,q})$ into $\frac{q}{2}$ edge-disjoint Hamilton cycles $C^1, C^2, \ldots, C^{q/2}$.
Let $$\mathcal{C}\colonequals \big\{C^1, C^2, \ldots, C^{\lceil(q-1)/2\rceil}\big\}.$$

For each $C^k\in \mathcal{C}$, we view $C^k$ as an oriented cycle whose orientation is given by \eqref{eq:Ck}.
The \emph{net voltage} $\alpha(C^k)$ of $C^k$ is the sum of voltages of its arcs, that is,
$$\alpha(C^k)\colonequals \alpha(v_{2k}u_1)+\alpha(u_1v_{2k+1})+\alpha(v_{2k+1}u_2)+\alpha(u_2v_{2k+2})+\cdots + \alpha(u_qv_{2k}).$$
Let $(C^k)^{\alpha}$ be the subgraph of $H_q$ lifted from $C^k$ with respect to the voltage assignment $\alpha$.
Since $C^k$ is a Hamilton cycle of $G_q$, $(C^k)^{\alpha}$ contains every vertex of $H_q$.
Lemma~\ref{le:cycle_lift} below provides the structure of each $(C^k)^{\alpha}$.
Before stating and proving the lemma, we first recall some facts on cosets.
For a subspace $L$ of the vector space $\mathbb{F}_q$ over $\mathbb{F}_p$, a \emph{coset} of $L$ is a set $x+L=\{x+t\colon\, t\in L\}$.
Two cosets are equal if their representatives differ by an element of $L$.
The cosets form a vector space $\mathbb{F}_q/L=\{x+L\colon\, x\in \mathbb{F}_q\}$, whose elements are the cosets rather than individual field elements.
If $L$ has $p$ elements, then $\mathbb{F}_q/L$ has $\frac{q}{p}$ elements and dimension $n-1$ over $\mathbb{F}_p$.
Moreover, the cosets of $L$ form a partition of $\mathbb{F}_q$.

\begin{lemma}\label{le:cycle_lift}
Let $C^k\in \mathcal{C}$ with net voltage $\alpha(C^k)=a\in \mathbb{F}_q$.
If $a= 0$, then $(C^k)^{\alpha}$ consists of $q$ vertex-disjoint cycles of length $2q$.
If $a\neq 0$, then $(C^k)^{\alpha}$ consists of $\frac{q}{p}$ vertex-disjoint cycles of length $2qp$.

Moreover, if $a\neq 0$, put $L=a\mathbb{F}_p$ and fix a vertex $w$ of $C^k$.
Then the components of $(C^k)^{\alpha}$ can be indexed by the $\frac{q}{p}$ cosets of $L$ in the following sense.
For any component of $(C^k)^{\alpha}$, there exists some $g\in \mathbb{F}_q$ such that this component meets $\phi(w)$ in exactly $\{(w,g), (w,g+a), \ldots, (w, g+(p-1)a)\}$,
and the second coordinates of these vertices form the coset $g+L$.
Thus there is a bijection between the components of $(C^k)^{\alpha}$ and the cosets of $L$.
\end{lemma}

\begin{proof}
Fix a vertex $w$ of $C^k$ arbitrarily.
For any $g\in \mathbb{F}_q$, we start at $(w,g)\in V(H_q)$ and follow the lift of repeated traversals of the oriented base cycle $C^k$.
After one traversal of the base cycle, we reach $(w,g+a)$, and after $j$ traversals, we reach $(w,g+ja)$.

If $a=0$, then the lifted walk closes after the first traversal.
Thus each component of $(C^k)^{\alpha}$ is a cycle of length $2q$.
Since there are $q$ different starting vertices $(w,g)$, $(C^k)^{\alpha}$ has $q$ different components.

If $a\neq 0$, then the smallest positive integer $j$ with $ja=0$ is $p$.
Thus each component of $(C^k)^{\alpha}$ is a cycle of length $2qp$, and its vertices in $\phi(w)$ are $(w,g), (w,g+a), \ldots, (w, g+(p-1)a)$.
Thus the second coordinates of these vertices form the coset $g+L$.
Moreover, two starting vertices $(w,g)$ and $(w,g')$ give the same component if and only if $g$ and $g'$ are in the same coset.
Since each coset has $p$ elements, $(C^k)^{\alpha}$ has $\frac{q}{p}$ different components, and each component can be indexed by one of the $\frac{q}{p}$ cosets of $L$.
\end{proof}

\subsection{Proof sketch}
\label{subsec:sketch}

Our starting point is the family $C^1, C^2, \ldots, C^{\lceil(q-1)/2\rceil}$ of edge-disjoint Hamilton cycles in the $K_{q,q}$ part of $G_q$ constructed in Section~\ref{subsec:cycles}.
By Lemma~\ref{le:cycle_lift}, their lifts $(C^1)^{\alpha}, \ldots, (C^{\lceil(q-1)/2\rceil})^{\alpha}$ are edge-disjoint spanning subgraphs of $H_q$, each consisting of some vertex-disjoint cycles.
We modify each lift into a Hamilton cycle of $H_q$ by deleting selected cross edges and adding edges within vertex fibres.
There are two difficulties: the components of each lift must be joined into one cycle, and the modifications for different lifts must remain edge-disjoint.

Suppose first that $a=\alpha(C^k)\neq 0$.
By Lemma~\ref{le:cycle_lift}, the components of $(C^k)^\alpha$ are indexed by the cosets of the one-dimensional subspace $L=a\mathbb{F}_p$.
When $q=p$, the lift is already a Hamilton cycle.
For $q=p^n$ with $n\geq 2$, Lemma~\ref{le:join} provides a construction using only the fibres above the two ends of a chosen base edge (see Figure~\ref{fig:cycle-joining}).
We delete two edges from each component and add one edge inside one vertex fibre, obtaining one spanning path from each component.
We then join these paths cyclically using edges inside the other vertex fibre.
To choose these edges, we use the fact that each difference set spans $\mathbb{F}_q$ to construct a cyclic ordering of the cosets from a Hamilton cycle in an additive Cayley graph.
The successive coset differences have representatives in the required difference set whose sum in $\mathbb{F}_q$ is zero.
This ensures that the joining construction closes (see Lemmas~\ref{le:Cayley}, \ref{le:z_i} and Corollary~\ref{cor:join}).

The case $\alpha(C^k)= 0$ requires different modifications because its lift has $q$ components, each containing only one vertex in every fibre.
We show in Lemma~\ref{le:alpha} that there is at most one base cycle with zero net voltage, occurring precisely when $q\equiv 3 \pmod{4}$ and $k=\frac{q+1}{4}$.
If $q=p$ is prime, we join its components successively using pairs of edges of difference $1$ within fibres of vertices (see Figure~\ref{fig:type2}).
If $q=p^n$ with $n>1$, we first group the components according to cosets of $\mathbb{F}_p$ and join the $p$ components in each group.
We then join the resulting cycles by the same joining procedure as in the proof of Lemma~\ref{le:join}, using two surviving lifted copies of one base edge, whose end-vertices in the fibre above one vertex have second coordinates differing by $2$.

To ensure edge-disjointness, for the case $\alpha(C^k)\neq 0$, we choose $e_k=u_{q-k}v_k$ and delete cross edges only above $e_k$.
The edges $e_k$ have pairwise distinct end-vertices as $k$ varies, so the added edges lie in distinct vertex fibres.
For the case $\alpha(C^k)= 0$, we choose the required edges from a matching while avoiding all fibres already used in other cases.
Moreover, since $C^1, \ldots, C^{\lceil(q-1)/2\rceil}$ are edge-disjoint in the base graph $G_q$, all the cross edges of $(C^1)^{\alpha}, \ldots, (C^{\lceil(q-1)/2\rceil})^{\alpha}$ are distinct.
Hence, the Hamilton cycles of $H_q$ just constructed are edge-disjoint.

\section{Joining components of a lifted Hamilton cycle}
\label{sec:join}

In Section~\ref{subsec:cycles}, we show that each Hamilton cycle $C^k$ of $G_q$ can be lifted to a subgraph $(C^k)^{\alpha}$ of $H_q$, where the components of $(C^k)^{\alpha}$ are cycles.
In this section, we show how to modify the components of $(C^k)^{\alpha}$ into a Hamilton cycle of $H_q$ for the case $\alpha(C^k)\neq 0$ and $n\geq 2$.
The idea is to delete two suitably chosen edges from each component and then use edges inside two fibres to join the resulting paths (see Lemma~\ref{le:join}).
Note that both $G_q$ and $H_q$ are undirected, and every ordinary edge has two opposite arcs.
In what follows, deleting either arc means deleting the corresponding undirected edge.

\begin{lemma}\label{le:join}
Let $q=p^n$ with $n\geq 2$, and let $C^k\in \mathcal{C}$ be a Hamilton cycle of $G_q$ with $\alpha(C^k)=a\neq 0$.
Let $L=a\mathbb{F}_p$, $m=\frac{q}{p}$, and $e\in E(C^k)$ with end-vertices $w$ and $t$.
Suppose that the following conditions hold:
\begin{itemize}
\item[{\rm (i)}] The difference set of $\phi(w)$ contains a nonzero element $s\in L$.
\item[{\rm (ii)}] The $m$ cosets of $L$ can be written in a cyclic order $A_0,A_1,\ldots,A_{m-1},A_m=A_0$, and there exist elements $z_0,z_1,\ldots,z_{m-1}$ in the difference set of $\phi(t)$ such that
    \begin{equation}\label{eq:splice}
    z_i+L=A_{i+1}-A_i \quad (0\leq i\leq m-1), \qquad \sum\nolimits_{i=0}^{m-1}z_i=0.
    \end{equation}
\end{itemize}
Then we can modify $(C^k)^{\alpha}$ into a Hamilton cycle of $H_q$ by deleting edges in the fibre above $e$ and adding edges inside the fibres above $w$ and $t$.
\end{lemma}

Before proving Lemma~\ref{le:join}, we first prove the following observation that will be used repeatedly.

\begin{observation}\label{obs:coset}
Let $q=p^n$ with $n\geq2$, let $L=a\mathbb{F}_p$ for some $a\neq 0$, and put $m=\frac{q}{p}$.
Suppose that the cosets of $L$ admit a cyclic ordering $A_0,\ldots,A_{m-1},A_m=A_0$, and that $z_0, z_1, \ldots, z_{m-1}\in\mathbb{F}_q$ satisfy~\eqref{eq:splice}.
For $s\in L$ and $x_0\in A_0$, define $x_{i+1}=x_i+z_i-s$ for $0\leq i\leq m-1$.
Then $x_i\in A_i$ for $0\leq i\leq m-1$, and $x_m=x_0$.
\end{observation}

\begin{proof}
We first prove $x_i\in A_i$ by induction on $i$.
The base case $x_0\in A_0$ is trivial.
Now assume that $x_i\in A_i$ and we shall show that $x_{i+1}\in A_{i+1}$.
Since $x_i\in A_i$, we have $x_{i+1}= x_i+z_i-s\in z_i-s+A_i$.
Since $s\in L$, $z_i+L=A_{i+1}-A_i$, and $A_i$ is a coset of $L$, we have $z_i-s+A_i=z_i+A_i=A_{i+1}$.
Hence, $x_{i+1}\in A_{i+1}$.
Furthermore, we have $x_m-x_0=\sum\nolimits_{i=0}^{m-1}(x_{i+1}-x_i)=\sum\nolimits_{i=0}^{m-1}(z_i-s)=\sum\nolimits_{i=0}^{m-1}z_i-ms$ by the definition of $x_i$.
Since $n\geq2$, the integer $m=\frac{q}{p}=p^{n-1}$ is divisible by $p$, so $ms=0$.
Combining with $\sum\nolimits_{i=0}^{m-1}z_i=0$, we have $x_m-x_0=\sum\nolimits_{i=0}^{m-1}z_i-ms=0$, so $x_m=x_0$.
\end{proof}

We now prove Lemma~\ref{le:join}.

\begin{figure}[htbp]
\begin{center}
\begingroup
\definecolor{JoinW}{RGB}{180,45,45}
\definecolor{JoinT}{RGB}{78, 143, 228}
\begin{tikzpicture}[scale=0.6,x=1cm,y=1cm,font=\small,
  retained/.style={draw=black,line width=0.85pt},
  removed/.style={draw=black,line width=0.85pt,dashed},
  addw/.style={draw=JoinW,line width=1.2pt},
  addt/.style={draw=JoinT,line width=1.2pt,preaction={draw=white,line width=3.6pt}},
  dot/.style={circle,draw=black,fill=black,inner sep=0pt,minimum size=3.8pt}]
\node at (4.2,2.55) {Vertex-disjoint cycles};
\node at (18.2,2.55) {One Hamilton cycle};
\begin{scope}[xshift=0cm]
\begin{scope}[yshift=0cm]
\coordinate (before-0-w0) at (0,0);
\coordinate (before-0-t0) at (2.8,0);
\coordinate (before-0-w1) at (5.4,0);
\coordinate (before-0-t1) at (8.4,0);
\draw[retained] (8.4,0) arc (-90:90:0.9) -- (0,1.8) arc (90:270:0.9);
\draw[retained] (2.8,0) -- (5.4,0);
\draw[removed] (0,0) -- (2.8,0);
\draw[removed] (5.4,0) -- (8.4,0);
\node at (4.2,1.3) {$(C^{k})^{\alpha}_{0}$};
\end{scope}
\begin{scope}[yshift=-3.7cm]
\coordinate (before-1-w0) at (0,0);
\coordinate (before-1-t0) at (2.8,0);
\coordinate (before-1-w1) at (5.4,0);
\coordinate (before-1-t1) at (8.4,0);
\draw[retained] (8.4,0) arc (-90:90:0.9) -- (0,1.8) arc (90:270:0.9);
\draw[retained] (2.8,0) -- (5.4,0);
\draw[removed] (0,0) -- (2.8,0);
\draw[removed] (5.4,0) -- (8.4,0);
\node at (4.2,1.3) {$(C^{k})^{\alpha}_{1}$};
\end{scope}
\begin{scope}[yshift=-8.3cm]
\coordinate (before-2-w0) at (0,0);
\coordinate (before-2-t0) at (2.8,0);
\coordinate (before-2-w1) at (5.4,0);
\coordinate (before-2-t1) at (8.4,0);
\draw[retained] (8.4,0) arc (-90:90:0.9) -- (0,1.8) arc (90:270:0.9);
\draw[retained] (2.8,0) -- (5.4,0);
\draw[removed] (0,0) -- (2.8,0);
\draw[removed] (5.4,0) -- (8.4,0);
\node at (4.2,1.3) {$(C^{k})^{\alpha}_{m-1}$};
\end{scope}
\node[font=\Large] at (4.2,-5.35) {$\vdots$};
\end{scope}
\begin{scope}[xshift=14cm]
\begin{scope}[yshift=0cm]
\coordinate (after-0-w0) at (0,0);
\coordinate (after-0-t0) at (2.8,0);
\coordinate (after-0-w1) at (5.4,0);
\coordinate (after-0-t1) at (8.4,0);
\draw[retained] (8.4,0) arc (-90:90:0.9) -- (0,1.8) arc (90:270:0.9);
\draw[retained] (2.8,0) -- (5.4,0);
\draw[addw] (0,0) .. controls (0.8,1.35) and (4.6,1.35) .. (5.4,0);
\node at (4.2,1.35) {$P(0)$};
\end{scope}
\begin{scope}[yshift=-3.7cm]
\coordinate (after-1-w0) at (0,0);
\coordinate (after-1-t0) at (2.8,0);
\coordinate (after-1-w1) at (5.4,0);
\coordinate (after-1-t1) at (8.4,0);
\draw[retained] (8.4,0) arc (-90:90:0.9) -- (0,1.8) arc (90:270:0.9);
\draw[retained] (2.8,0) -- (5.4,0);
\draw[addw] (0,0) .. controls (0.8,1.35) and (4.6,1.35) .. (5.4,0);
\node at (4.2,1.35) {$P(1)$};
\end{scope}
\begin{scope}[yshift=-8.3cm]
\coordinate (after-2-w0) at (0,0);
\coordinate (after-2-t0) at (2.8,0);
\coordinate (after-2-w1) at (5.4,0);
\coordinate (after-2-t1) at (8.4,0);
\draw[retained] (8.4,0) arc (-90:90:0.9) -- (0,1.8) arc (90:270:0.9);
\draw[retained] (2.8,0) -- (5.4,0);
\draw[addw] (0,0) .. controls (0.8,1.35) and (4.6,1.35) .. (5.4,0);
\node at (5.2,1.35) {$P(m-1)$};
\end{scope}
\node[font=\Large] at (4.2,-5.35) {$\vdots$};
\end{scope}
\draw[->,line width=1.1pt] (10,-3.3) -- (12.5,-3.3);
\draw[addt] (after-0-t0) -- (after-1-t1);
\draw[addt] (after-1-t0) -- (18.8,-5.342857);
\node[rotate=-39.4,text=JoinT,fill=white,inner sep=1pt] at (19.4,-5.835714) {$\cdots$};
\draw[addt] (20.0,-6.328571) -- (after-2-t1);
\draw[addt] (after-2-t0) -- (after-0-t1);
\node[dot] at (before-0-w0) {};
\node[dot] at (before-0-t0) {};
\node[dot] at (before-0-w1) {};
\node[dot] at (before-0-t1) {};
\node[dot] at (before-1-w0) {};
\node[dot] at (before-1-t0) {};
\node[dot] at (before-1-w1) {};
\node[dot] at (before-1-t1) {};
\node[dot] at (before-2-w0) {};
\node[dot] at (before-2-t0) {};
\node[dot] at (before-2-w1) {};
\node[dot] at (before-2-t1) {};
\node[dot] at (after-0-w0) {};
\node[dot] at (after-0-t0) {};
\node[dot] at (after-0-w1) {};
\node[dot] at (after-0-t1) {};
\node[dot] at (after-1-w0) {};
\node[dot] at (after-1-t0) {};
\node[dot] at (after-1-w1) {};
\node[dot] at (after-1-t1) {};
\node[dot] at (after-2-w0) {};
\node[dot] at (after-2-t0) {};
\node[dot] at (after-2-w1) {};
\node[dot] at (after-2-t1) {};
\begin{scope}[xshift=0cm]
\node at (0.2,0.50) {$(w,x_0)$};
\node at (2.7,-0.50) {$(t,x_0+\beta)$};
\node at (5.2,0.50) {$(w,x_0+s)$};
\node at (8.50,-0.50) {$(t,x_0+s+\beta)$};
\end{scope}
\begin{scope}[yshift=-3.7cm]
\node at (0.2,0.50) {$(w,x_1)$};
\node at (2.7,-0.50) {$(t,x_1+\beta)$};
\node at (5.2,0.50) {$(w,x_1+s)$};
\node at (8.50,-0.50) {$(t,x_1+s+\beta)$};
\end{scope}
\begin{scope}[yshift=-8.3cm]
\node at (0.6,0.50) {$(w,x_{m-1})$};
\node at (2.7,-0.50) {$(t,x_{m-1}+\beta)$};
\node at (5.2,0.50) {$(w,x_{m-1}+s)$};
\node at (8.50,-0.50) {$(t,x_{m-1}+s+\beta)$};
\end{scope}
\draw[removed] (-0.7,-10) -- (0.15,-10);
\node[anchor=west] at (0.3,-10) {deleted edges};
\draw[addw] (6.2,-10) -- (7.05,-10);
\node[anchor=west] at (7.2,-10) {added edges within $\phi(w)$};
\draw[addt] (15.8,-10) -- (16.65,-10);
\node[anchor=west] at (16.8,-10) {added edges within $\phi(t)$};
\end{tikzpicture}
\endgroup
\caption{Joining the cycles $(C^k)^{\alpha}_0, (C^k)^{\alpha}_1, \ldots, (C^k)^{\alpha}_{m-1}$ in the proof of Lemma~\ref{le:join}.}
\label{fig:cycle-joining}
\end{center}
\vspace{-0.5cm}
\end{figure}
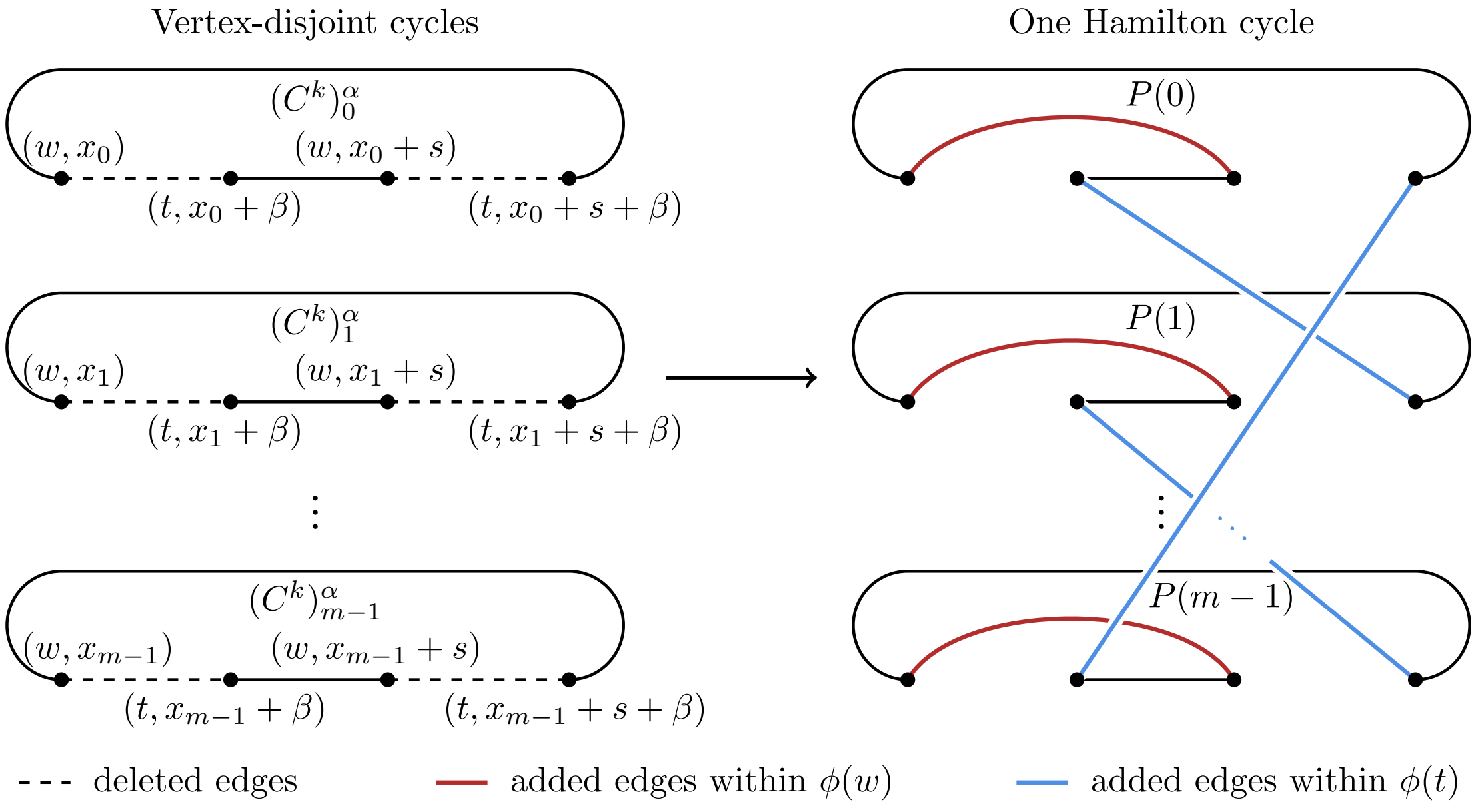

\begin{proof}[Proof of Lemma~\ref{le:join}]
By Lemma~\ref{le:cycle_lift}, $(C^k)^{\alpha}$ consists of $m$ vertex-disjoint cycles that can be indexed by the $m$ cosets $A_0(=A_m), A_1, \ldots, A_{m-1}$ of $L$.
For each $i\in \{0, 1, \ldots, m-1\}$, we write $(C^k)^{\alpha}_i$ for the cycle in $(C^k)^{\alpha}$ such that the second coordinates of the vertices in $V\big((C^k)^{\alpha}_i\big)\cap \phi(w)$ form the coset $A_i$.
Take an $x_0\in A_0$ and define successively $x_{i+1}\colonequals x_i+z_i-s$ for $0\leq i\leq m-1$.
By Observation~\ref{obs:coset}, we have $x_i\in A_i$ for $0\leq i\leq m-1$, and $x_m=x_0$.

We now join the cycles $(C^k)^{\alpha}_0, (C^k)^{\alpha}_1, \ldots, (C^k)^{\alpha}_{m-1}$ into a Hamilton cycle of $H_q$ (see Figure~\ref{fig:cycle-joining} for an illustration).
For each $i\in \{0, 1, \ldots, m-1\}$, we have $x_i+s\in A_i$ since $x_i\in A_i$ and $s\in L$.
Thus $(C^k)^{\alpha}_i$ contains two distinct edges in $\phi(e)$ whose end-vertices in $\phi(w)$ are $(w,x_i)$ and $(w,x_i+s)$.
Let $\beta$ denote the voltage of the arc from $w$ to $t$.
Then the other end-vertices of these two edges are $(t, x_i+\beta)$ and $(t, x_i+s+\beta)$, respectively.
The two chosen edges lie above the same base edge $e$, so they have the same direction in the orientation inherited from $C^k$.
We delete these two edges from $(C^k)^{\alpha}_i$.
Then we obtain two paths, one with end-vertices $(t, x_i+\beta)$ and $(w, x_i+s)$, and the other one with end-vertices $(t, x_i+s+\beta)$ and $(w, x_i)$.
Now we add the edge between $(w, x_i)$ and $(w, x_i+s)$, which exists by condition~(i).
This joins the two paths into one path $P(i)$ with end-vertices $(t, x_i+\beta)$ and $(t, x_i+s+\beta)$, and $P(i)$ contains every vertex of $(C^k)^{\alpha}_i$.

It remains to join the $m$ paths $P(0), P(1), \ldots, P(m-1)$ into a Hamilton cycle of $H_q$.
For every $i$ (with indices modulo $m$), we add an edge from $(t,x_i+\beta)$ to $(t,x_{i+1}+s+\beta)$.
Such an edge exists since the difference between its second coordinates is $x_{i+1}+s-x_i=z_i$, which is contained in the difference set of $\phi(t)$ by condition~(ii).
This joins the paths in the cyclic order $P(0), P(1), \ldots, P(m-1), P(0)$.
The resulting graph is a cycle containing every vertex of $(C^k)^{\alpha}$ (and thus every vertex of $H_q$), so it is a Hamilton cycle of $H_q$.
This completes the proof of Lemma~\ref{le:join}.
\end{proof}

Once condition~(i) of Lemma~\ref{le:join} holds, it remains to find a cyclic ordering of the cosets of $L$ and elements $z_i$ satisfying condition~(ii).
Lemma~\ref{le:z_i} supplies these under a spanning hypothesis on the difference set of $\phi(t)$.
In order to prove Lemma~\ref{le:z_i}, we first prove the following result on Cayley graphs.
For a vector space $W$ and a subset $B\subseteq W$ with $B=-B$, the \emph{Cayley graph} $\operatorname{Cay}(W,B)$ has vertex set $W$, and two distinct vertices are adjacent if and only if their difference belongs to $B$.

\begin{lemma}\label{le:Cayley}
Let $W$ be a $d$-dimensional vector space over $\mathbb{F}_p$ with basis $b_1,\ldots,b_d$.
If $p$ is odd and $d\geq 1$, or if $p=2$ and $d\geq 2$, then $\operatorname{Cay}(W,\{\pm b_1,\ldots,\pm b_d\})$ contains a Hamilton cycle.
\end{lemma}

\begin{proof}
We prove the assertion by induction on $d$.
For odd $p$ and $d=1$, the sequence $0,b_1,2b_1,\ldots,(p-1)b_1,0$ is a Hamilton cycle.
For $p=2$ and $d=2$, the sequence $0,b_1,b_1+b_2,b_2,0$ is a Hamilton cycle.
Now assume that the result holds for $d-1$, and we shall prove it for $d$.

Let $W'=\operatorname{span}_{\mathbb{F}_p}\{b_1,\ldots,b_{d-1}\}\subseteq W.$
Then $|W'|=p^{d-1}$, and $W', W'+b_d, \ldots, W'+(p-1)b_d$ form a partition of $W$.
By the induction hypothesis, $\operatorname{Cay}(W',\{\pm b_1,\ldots,\pm b_{d-1}\})$ contains a Hamilton cycle $x_0x_1\cdots x_{M-1}x_{M}$, where $x_{M}=x_0$ and $M=p^{d-1}$.
For each $j\in \{1, 2, \ldots, M\}$, let $h_j=x_j-x_{j-1}$ be the difference between two successive vertices.
Then $h_1, h_2, \ldots, h_M\in \{\pm b_1,\ldots,\pm b_{d-1}\}$ and $\sum\nolimits_{j=1}^{M}h_j=0$.
Note that by deleting the edge $x_{M-1}x_{M}$, we obtain a Hamilton path $x_0x_1\cdots x_{M-1}$ in $\operatorname{Cay}(W',\{\pm b_1,\ldots,\pm b_{d-1}\})$ with successive differences $h_1, h_2, \ldots, h_{M-1}$.

Put $t=b_d-h_M$.
For each $i\in\{0,1,\ldots,p-1\}$, consider the path
$$Q(i):\quad x_0+it,\ x_1+it,\ \ldots,\ x_{M-1}+it$$
in $\operatorname{Cay}(W,\{\pm b_1,\ldots,\pm b_d\})$.
The successive differences along $Q(i)$ are $h_1,\ldots,h_{M-1}$, so consecutive vertices are adjacent in $\operatorname{Cay}(W,\{\pm b_1,\ldots,\pm b_d\})$.
Since $h_M\in W'$, we have $W'+it=W'+i(b_d-h_M)=W'+ib_d$.
Thus $Q(i)$ visits every vertex of $W'+ib_d$ exactly once.
Hence, $Q(i)$ is a Hamilton path in the subgraph of $\operatorname{Cay}(W,\{\pm b_1,\ldots,\pm b_d\})$ induced by $W'+ib_d$.
We claim that the paths $Q(0),Q(1),\ldots,Q(p-1)$ can be joined successively into a Hamilton cycle of $\operatorname{Cay}(W,\{\pm b_1,\ldots,\pm b_d\})$.
Indeed, for each $i\in \{0,1,\ldots,p-2\}$, the difference between the final vertex of $Q(i)$ and the initial vertex of $Q(i+1)$ is
$$(x_0+(i+1)t)-(x_{M-1}+it) =x_0-x_{M-1}+t =h_M+b_d-h_M=b_d.$$
Hence, these two vertices are adjacent.
Moreover, since $pt=0$, we have
$$x_0-(x_{M-1}+(p-1)t) =x_0-x_{M-1}+t =h_M+t =b_d.$$
Thus the final vertex of $Q(p-1)$ is adjacent to the initial vertex of $Q(0)$.
Recall that $W',W'+b_d,\ldots,W'+(p-1)b_d$ form a partition of $W$.
Hence, the paths $Q(0),Q(1),\ldots,Q(p-1)$ are pairwise vertex-disjoint and together contain every vertex of $W$.
By joining consecutive paths using the edges described above, we obtain a Hamilton cycle in $\operatorname{Cay}(W,\{\pm b_1,\ldots,\pm b_d\})$.
This completes the proof.
\end{proof}

\begin{lemma}\label{le:z_i}
Let $q=p^n$ with $n\geq 2$.
Let $L=a\mathbb{F}_p$ for some $a\in \mathbb{F}_q\setminus\{0\}$, put $m=\frac{q}{p}$, and let $S\subseteq \mathbb{F}_q^{\times}$ with $S=-S$.
Suppose that the images of the elements of $S$ span $\mathbb{F}_q/L$, i.e., $\operatorname{span}_{\mathbb{F}_p}\{s+L\colon\, s\in S\}=\mathbb{F}_q/L=\{x+L\colon\, x\in \mathbb{F}_q\}$.
Then the cosets of $L$ admit a cyclic ordering $A_0,A_1,\ldots,A_{m-1},A_m=A_0$ such that there exist $z_0,z_1, \ldots,z_{m-1}\in S$ satisfying
\begin{equation*}\label{eq:quotient}
z_i+L=A_{i+1}-A_i \quad (0\leq i\leq m-1), \qquad \sum\nolimits_{i=0}^{m-1}z_i=0.
\end{equation*}
\end{lemma}

\begin{proof}
We first consider the case $q=2^2=4$.
Since the images of $S$ span $\mathbb{F}_4/L$, there exists an element $b\in S\setminus L$.
Put $A_0=L$, $A_1=b+L$, $A_2=A_0$, and $z_0=z_1=b$.
Then $A_0=A_2$ and $A_1$ are the two distinct cosets of $L$.
Since the characteristic of $\mathbb{F}_4$ is 2, we have $z_i+L=A_{i+1}-A_i$ for $i\in \{0,1\}$, and $z_0+z_1=b+b=0$.
Thus the conclusion holds.

Now let $q=p^n$ with $n\geq 2$ when $p$ is odd, and $n\geq 3$ when $p=2$.
Since $\mathbb{F}_q/L$ has dimension $n-1$ and $\{s+L\colon\, s\in S\}$ spans $\mathbb{F}_q/L$, we may choose $b_1, \ldots, b_{n-1}\in S$ such that the cosets $b_1+L, \ldots, b_{n-1}+L$ form a basis of $\mathbb{F}_q/L$.
Let $T=\operatorname{span}_{\mathbb{F}_p}\{b_1, \ldots, b_{n-1}\}$.

\begin{claim}\label{cl:z_i}
The vectors $b_1,\ldots,b_{n-1}$ are linearly independent, and $T$ contains exactly one element of each coset of $L$.
\end{claim}

\begin{proof}
We first show that the vectors $b_1,\ldots,b_{n-1}$ are linearly independent.
Indeed, if $c_1b_1+\cdots+c_{n-1}b_{n-1}=0$ with $c_1, \ldots, c_{n-1}\in \mathbb{F}_p$, then $c_1(b_1+L)+\cdots+c_{n-1}(b_{n-1}+L)=L$, which is the zero element of $\mathbb{F}_q/L$.
Thus $c_1=\cdots=c_{n-1}=0$, so the vectors $b_1,\ldots,b_{n-1}$ are linearly independent.
In particular, $T$ has $p^{n-1}=m$ elements.

Moreover, we claim that the $m$ elements of $T$ belong to $m$ different cosets.
To see this, assume that two elements $x,y\in T$ belong to the same coset of $L$, i.e., $x+L=y+L$.
Then $x-y\in L$.
Write $x-y= \sum\nolimits_{j=1}^{n-1}c'_jb_j$ for some $c'_1, \ldots, c'_{n-1}\in \mathbb{F}_p$.
Then $x-y+L=\sum\nolimits_{j=1}^{n-1}c'_jb_j+L= \sum\nolimits_{j=1}^{n-1}c'_j(b_j+L)$.
Since $x-y\in L$, we have $x-y+L=L$, the zero element of $\mathbb{F}_q/L$.
Thus $c'_1=\cdots=c'_{n-1}=0$, so $x=y$.
Hence, the $m$ elements of $T$ belong to $m$ different cosets.
Since there are exactly $m$ cosets of $L$ in $\mathbb{F}_q$, $T$ contains exactly one element of each coset of $L$.
\end{proof}

By Claim~\ref{cl:z_i} and Lemma~\ref{le:Cayley}, $\operatorname{Cay}(T,\{\pm b_1,\ldots,\pm b_{n-1}\})$ contains a Hamilton cycle.
Write it as $x_0x_1\cdots x_{m-1}x_{m}$, where $x_{m}=x_0$.
Let $A_i=x_i+L$ for $0\leq i\leq m$, and let $z_i=x_{i+1}-x_i$ for $0\leq i\leq m-1$.
Then $z_0, z_1, \ldots, z_{m-1}\in \{\pm b_1,\ldots,\pm b_{n-1}\}\subseteq S$, $z_i+L=A_{i+1}-A_i$, and $A_0,A_1,\ldots,A_{m-1}$ list all cosets of $L$ by Claim~\ref{cl:z_i}.
Moreover, $\sum\nolimits_{i=0}^{m-1}z_i=\sum\nolimits_{i=0}^{m-1}(x_{i+1}-x_i)=x_m-x_0=0$ in $\mathbb{F}_q$.
This completes the proof of Lemma~\ref{le:z_i}.
\end{proof}

Applying Lemmas~\ref{le:difference}, \ref{le:join} and \ref{le:z_i}, we obtain the following result.

\begin{corollary}\label{cor:join}
Let $q=p^n$ with $n\geq 2$.
Let $C^k\in \mathcal{C}$ with $\alpha(C^k)\neq 0$, and $e\in E(C^k)$ with end-vertices $u\in U$ and $v\in V$.
Then we can modify $(C^k)^{\alpha}$ into a Hamilton cycle of $H_q$ by deleting edges in the fibre above $e$ and adding edges inside the fibres above $u$ and $v$.
\end{corollary}

\begin{proof}
Let $a=\alpha(C^k)$ and $L=a\mathbb{F}_p$.
Recall that the difference sets $S_U$ and $S_V$ satisfy $S_U=-S_U$ and $S_V=-S_V$.
By Lemma~\ref{le:difference}, we have $S_U\subseteq \mathbb{F}_q^{\times}$, $S_V\subseteq \mathbb{F}_q^{\times}$, and each of $S_U$ and $S_V$ spans $\mathbb{F}_q$ as a vector space over $\mathbb{F}_p$.
Then the images of each of $S_U$ and $S_V$ span $\mathbb{F}_q/L$.
Hence, Lemma~\ref{le:z_i} applies to both $S_U$ and $S_V$.

Moreover, we have $S_U\cup S_V=\mathbb{F}_q^{\times}$ by Lemma~\ref{le:difference}.
Thus the nonzero element $a$ belongs to at least one of $S_U$ and $S_V$.
If $a\in S_U$, take $w=u$, $t=v$, and $s=a$; otherwise take $w=v$, $t=u$, and $s=a$.
Now condition~(i) of Lemma~\ref{le:join} holds.
Applying Lemma~\ref{le:z_i} to the difference set of the fibre above $t$, we have that condition~(ii) of Lemma~\ref{le:join} holds.
Then the result follows from Lemma~\ref{le:join}.
\end{proof}

\section{Edge-disjoint Hamilton cycles in $H_q$}
\label{sec:proof}

In this section, we present our proof of Theorem~\ref{thm:q}.
We first determine which of the cycles $C^k$ in $\mathcal{C}$ have zero net voltage, extending \cite[Claims~3.1 and 3.2]{LiWZ}.

\begin{lemma}\label{le:alpha}
For every prime power $q\geq 3$ and every integer $k$ with $1\leq k\leq \big\lceil\frac{q-1}{2}\big\rceil$, $\alpha(C^k)=0$ if and only if $q=4k-1$.
\end{lemma}

\begin{proof}
For each $C^k\in \mathcal{C}$, consider the orientation of $C^k$ given by \eqref{eq:Ck}:
\begin{equation*}\label{eq:Ck+}
v_{2k}\to u_1\to v_{2k+1}\to u_2\to v_{2k+2}\to \cdots \to u_q \to v_{2k},
\end{equation*}
where the indices in $u_i$ and $v_i$ are considered modulo $q$ and represented in $\{1, 2, \ldots, q\}$.
Recall that an arc from $v_j$ to $u_i$ has voltage $-\xi^i\xi^j=-\xi^{i+j}$ if both indices are different from $q$, whereas an arc from $u_i$ to $v_j$ has voltage $\xi^i\xi^j=\xi^{i+j}$.
If either index equals $q$, then the voltage is zero.

We first consider the case $q=3$.
In this case, we may take $\xi=2=-1$.
There is only one cycle $C^1=v_2u_1v_3u_2v_1u_3v_2$ in $\mathcal{C}$.
The voltages of its arcs are $-2, 0, 0, 2, 0, 0$, so $\alpha(C^1)=0$.

We next consider the case $k=\frac{q}{2}$.
In this case, $q\geq 4$ is even, and $C^{q/2}$ uses the matchings $M_q$ and $M_1$.
For $M_1$, since $\xi^2\neq 1$ and $\xi^{q-1}=1$, the sum of voltages of arcs is
$$\sum\nolimits_{i=1}^{q}\alpha(u_iv_i)=\sum\nolimits_{i=1}^{q-1}\xi^{i}\xi^{i}=\sum\nolimits_{i=1}^{q-1}\xi^{2i}=\frac{\xi^2(1-\xi^{2(q-1)})}{1-\xi^2}=\frac{\xi^2(1-1^2)}{1-\xi^2}=0.$$
For $M_q$, since the characteristic is 2, we have $-\xi^{2i-1}=\xi^{2i-1}$ and $-\xi=\xi$, so the sum of voltages of arcs is
$$\sum\nolimits_{i=1}^{q}\alpha(v_{i+q-1}u_i)=\sum\nolimits_{i=2}^{q-1}-\xi^{i-1}\xi^{i}=\sum\nolimits_{i=2}^{q-1}\xi^{2i-1}
=\xi+\sum\nolimits_{i=1}^{q-1}\xi^{2i-1}=\xi+\frac{\xi(1-\xi^{2(q-1)})}{1-\xi^2}=\xi.$$
Hence, $\alpha(C^{q/2})=\xi\ne 0$.

It remains to consider the case $q>3$ and $1\leq k< \frac{q}{2}$.
For each $1\leq i\leq q$, let $e_{2i-1}$ be the arc from $v_{i+2k-1}$ to $u_i$, and let $e_{2i}$ be the arc from $u_i$ to $v_{i+2k}$.
First, consider the arcs $e_{2i-1}$.
If $1\leq i\leq q-2k$, then $i+2k-1\leq q-1$, and thus $\alpha(e_{2i-1}) =-\xi^i\xi^{i+2k-1} =-\xi^{2i+2k-1}.$
If $i=q-2k+1$, then $v_{i+2k-1}=v_q$, so $\alpha(e_{2i-1})=0$.
If $q-2k+2\leq i\leq q-1$, then the index $i+2k-1$ is reduced modulo $q$, so $\alpha(e_{2i-1}) =-\xi^i\xi^{i+2k-1-q} =-\xi^{2i+2k-1-q}.$
Finally, $\alpha(e_{2q-1})=0$ since $e_{2q-1}$ is incident with $u_q$.
Next, consider the arcs $e_{2i}$.
If $1\leq i\leq q-2k-1$, then $i+2k\leq q-1$, and thus $\alpha(e_{2i})=\xi^i\xi^{i+2k}=\xi^{2i+2k}.$
If $i=q-2k$, then $v_{i+2k}=v_q$, so $\alpha(e_{2i})=0$.
If $q-2k+1\leq i\leq q-1$, then the index $i+2k$ is reduced modulo $q$, so $\alpha(e_{2i})=\xi^i\xi^{i+2k-q}=\xi^{2i+2k-q}.$
Finally, $\alpha(e_{2q})=0$ since $e_{2q}$ is incident with $u_q$.

The voltages of $e_1, e_2, \ldots, e_{2q}$ are
$$-\xi^{2k+1}, \xi^{2k+2}, \ldots, -\xi^{2q-2k-1}, 0, 0, \xi^{q-2k+2}, -\xi^{q-2k+3}, \ldots, \xi^{q+2k-2}, 0, 0.$$
Note that $\xi^j=-(-\xi)^j$ for odd $j$, and $\xi^j=(-\xi)^j$ for even $j$.
Thus for odd $q$, the voltages of $e_1, e_2, \ldots, e_{2q}$ equal
$$(-\xi)^{2k+1}, (-\xi)^{2k+2}, \ldots, (-\xi)^{2q-2k-1}, 0, 0, -(-\xi)^{q-2k+2}, -(-\xi)^{q-2k+3}, \ldots, -(-\xi)^{q+2k-2}, 0, 0.$$
This also holds for even $q$, since $-1=1$ in $\mathbb{F}_{2^n}.$
Therefore, we have
\begin{equation}\label{eq:voltage-sums}
\alpha(C^k) =\sum\nolimits_{j=2k+1}^{2q-2k-1}(-\xi)^j -\sum\nolimits_{j=q-2k+2}^{q+2k-2}(-\xi)^j.
\end{equation}
Since $q>3$ and $\xi$ is a primitive element of $\mathbb{F}_q$, we have $\xi \neq -1$, and hence $-\xi \neq 1$.
Applying the geometric-sum formula to \eqref{eq:voltage-sums}, we obtain
\begin{align*}
\alpha(C^k) = &~\frac{(-\xi)^{2q-2k}-(-\xi)^{2k+1}}{(-\xi)-1} -\frac{(-\xi)^{q+2k-1}-(-\xi)^{q-2k+2}}{(-\xi)-1} \\
= &~\frac{(-\xi)^{2q-2k}-(-\xi)^{2k+1} -(-\xi)^{q+2k-1}+(-\xi)^{q-2k+2}}{(-\xi)-1}.
\end{align*}
Moreover, we have $(-\xi)^{q-1}=\xi^{q-1}=1$.
Indeed, for odd $q$, the exponent $q-1$ is even, whereas for even $q$ we have $-\xi=\xi$.
Hence, $(-\xi)^{q+2k-1}=(-\xi)^{2k}$ and $(-\xi)^{q-2k+2}=(-\xi)^{2q-2k+1}.$
It follows that
\begin{align*}
\alpha(C^k)= &~\frac{(-\xi)^{2q-2k}-(-\xi)^{2k+1} -(-\xi)^{2k}+(-\xi)^{2q-2k+1}}{(-\xi)-1} = \frac{(1+(-\xi))\left((-\xi)^{2q-2k}-(-\xi)^{2k}\right)}{(-\xi)-1} \\
= &~\frac{(1+(-\xi))(-\xi)^{2k}\left((-\xi)^{2q-4k}-1\right)}{(-\xi)-1} = \frac{(1+(-\xi))(-\xi)^{2k}\left(\xi^{2q-4k}-1\right)}{(-\xi)-1}.
\end{align*}
Since $-\xi\notin \{-1, 0, 1\}$, this implies that $\alpha(C^k)=0$ if and only if $\xi^{2q-4k}=1$.
Recall that $\xi$ has multiplicative order $q-1$.
Hence, $\xi^{2q-4k}=1$ if and only if $q-1$ divides $2q-4k$.
Moreover, since $1\leq k< \frac{q}{2}$, we have $0< 2q-4k\leq 2q-4 <2(q-1)$.
Thus the only possible multiple of $q-1$ in this range is $q-1$.
Thus $\alpha(C^k)=0$ if and only if $2q-4k=q-1$, equivalently $q=4k-1$.
This completes the proof of Lemma~\ref{le:alpha}.
\end{proof}

Now we have all the ingredients to present our proof of Theorem~\ref{thm:q}.

\begin{proof}[Proof of Theorem~\ref{thm:q}]
We start from the edge-disjoint Hamilton cycles $C^1, C^2, \ldots, C^{\lceil(q-1)/2\rceil}$ of the base graph $G_q$ constructed in Section~\ref{subsec:cycles}.
We shall modify their lifts $(C^1)^{\alpha}, \ldots, (C^{\lceil(q-1)/2\rceil})^{\alpha}$ into edge-disjoint Hamilton cycles of $H_q$ one by one.
By Lemma~\ref{le:alpha}, $\alpha(C^k)=0$ if and only if $q=4k-1$.
Hence, $\alpha(C^k)=0$ occurs only when $q=4\ell-1$ (i.e., $q$ is odd and $q\equiv 3 \pmod{4}$) and $k=\ell$.

We first consider the case that $q$ is even.
In this case, $\alpha(C^k)\neq 0$ for every $k\in \{1, 2, \ldots, \frac{q}{2}\}$.
For each $k\in \{1, 2, \ldots, \frac{q}{2}\}$, let $e_k=u_{q-k}v_k$.
Then $e_k\in M_{2k+1}\subseteq E(C^k)$ since $v_{(q-k)+(2k+1)-1}=v_{q+k}=v_k$, where $M_{q+1}=M_1$.
By Corollary~\ref{cor:join}, we can modify $(C^k)^{\alpha}$ into a Hamilton cycle of $H_q$ by deleting edges in the fibre above $e_k$ and adding edges inside the fibres above $u_{q-k}$ and $v_{k}$.
As $k$ varies from $1$ to $\frac{q}{2}$, the vertices $u_{q-k}$ are distinct, and the vertices $v_k$ are distinct.
Hence, for different modifications, we add edges within distinct vertex fibres, so the added edges are disjoint.
Moreover, the cross edges are disjoint since the base cycles are edge-disjoint.
Therefore, we obtain $\frac{q}{2}$ pairwise edge-disjoint Hamilton cycles of $H_q$.
This completes the proof for even $q$.

In the remainder of the proof, we assume that $q$ is odd.
Note that if $q=p^n \equiv 3 \pmod{4}$, then $p\equiv 3 \pmod{4}$ and $n$ is odd.
Hence, it suffices to consider the following three types of cycles.

\smallskip\noindent\textbf{Type~1:} $\alpha(C^k)\neq 0$.\smallskip

By Lemma~\ref{le:cycle_lift}, $(C^k)^{\alpha}$ consists of $\frac{q}{p}$ vertex-disjoint cycles of length $2qp$.
If $n=1$, then $q=p$, and thus $(C^k)^{\alpha}$ itself is a Hamilton cycle of $H_q$.
Now assume $n\geq 2$.
In this case, the modification is analogous to that for even $q$.
For each $k\in \{1, 2, \ldots, \frac{q-1}{2}\}$, we choose the edge $e_k=u_{q-k}v_k\in M_{2k+1}\subseteq E(C^k)$.
By Corollary~\ref{cor:join}, we can modify $(C^k)^{\alpha}$ into a Hamilton cycle of $H_q$ by deleting arcs in the fibre above $e_k$ and adding edges inside the fibres above $u_{q-k}$ and $v_k$.
Note that for different $k\in \{1, 2, \ldots, \frac{q-1}{2}\}$ with $\alpha(C^k)\neq 0$, the vertices $u_{q-k}$ are distinct, and the vertices $v_k$ are distinct.
Hence, for different modifications, we add edges within distinct vertex fibres.
Moreover, the cross edges are also disjoint, since the base cycles $C^k$ are edge-disjoint.
Therefore, the Hamilton cycles constructed in this case are edge-disjoint.

\smallskip\noindent\textbf{Type~2:} $\alpha(C^{\ell})=0$ and $n=1$.\smallskip

In this case, we have $q=p$.
By Lemma~\ref{le:cycle_lift}, $(C^{\ell})^{\alpha}$ consists of $p$ vertex-disjoint cycles $C(0), C(1), \ldots, C(p-1)$ of length $2p$, each of the form
\begin{align*}\label{eq:Cl}
C(j):~&~(v_{2\ell},j)\to (u_1,j+\alpha(v_{2\ell}u_1))\to (v_{2\ell+1},j+\alpha(v_{2\ell}u_1)+\alpha(u_1v_{2\ell+1}))\to \\
~&~(u_2,j+\alpha(v_{2\ell}u_1)+\alpha(u_1v_{2\ell+1})+\alpha(v_{2\ell+1}u_2))\to \\
~&~(v_{2\ell+2},j+\alpha(v_{2\ell}u_1)+\alpha(u_1v_{2\ell+1})+\alpha(v_{2\ell+1}u_2)+\alpha(u_2v_{2\ell+2}))\to \cdots \to \\
~&~ (u_p,j+\alpha(v_{2\ell}u_1)+\alpha(u_1v_{2\ell+1})+\alpha(v_{2\ell+1}u_2)+\alpha(u_2v_{2\ell+2})+\dots+\alpha(v_{2\ell-1}u_p)) \to (v_{2\ell},j)
\end{align*}
with $j\in \mathbb{F}_p$.
Note that for each $1\leq j\leq p-1$, the second coordinates of vertices of $C(j-1)$ and $C(j)$ in the same vertex fibre have difference 1.

We join these cycles successively (see Figure~\ref{fig:type2} for an illustration).
For each $1\leq j\leq p-1$, we choose the edge $f_j=u_jv_{j+2\ell}\in M_{2\ell+1}$.
These $p-1$ edges are distinct.
Start with $C(0)$.
In step $j$, delete the lifted copy of $f_j$ belonging to the original component $C(j-1)$ from the cycle already constructed, and delete the lifted copy of $f_j$ in $C(j)$.
These edges are still present since the base edges $f_1, \ldots, f_{p-1}$ are distinct.
We obtain two paths with end-vertices in $\phi(u_j)$ and $\phi(v_{j+2\ell})$.
Join the two end-vertices in $\phi(u_j)$ and the two end-vertices in $\phi(v_{j+2\ell})$, respectively.
Their second coordinates differ by $1$, so both added edges exist since $1\in S_U\cap S_V$ by Lemma~\ref{le:difference} (note that we have $\delta=-1$ since $q\equiv 3 \pmod{4}$ in this case).
Then we obtain a single cycle on $V\big(C(0)\big)\cup \cdots \cup V\big(C(j)\big)$.
After $p-1$ steps, this is a Hamilton cycle of $H_p$.
The new edges are all distinct since $f_1, f_2, \ldots, f_{p-1}$ belong to a perfect matching in $G_p$.
They cannot overlap any other Hamilton cycles constructed for Type~1, since those cycles use only cross edges when $n=1$.

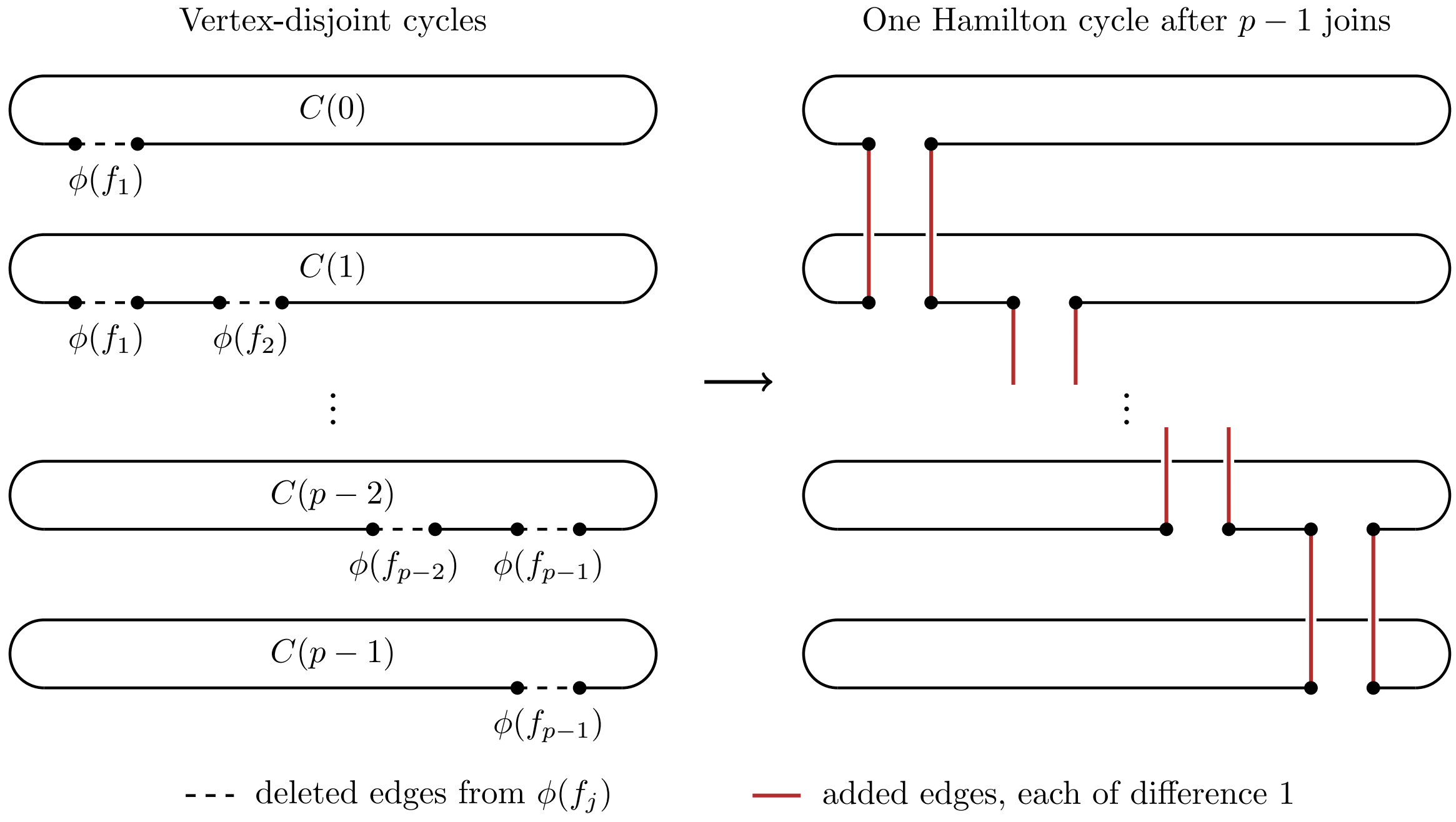
\begin{figure}[htbp]
\begin{center}
\begingroup
\definecolor{TypeTwoJoin}{RGB}{180,45,45}
\begin{tikzpicture}[scale=0.6,x=1cm,y=1cm,font=\small,
  retained/.style={draw=black,line width=0.85pt},
  removed/.style={draw=black,line width=0.85pt,dashed},
  joined/.style={draw=TypeTwoJoin,line width=1.15pt,preaction={draw=white,line width=3.3pt}},
  dot/.style={circle,fill=black,draw=black,inner sep=0pt,minimum size=3.7pt}]
\node at (5.1,2.15) {Vertex-disjoint cycles};
\node at (19.1,2.15) {One Hamilton cycle after $p-1$ joins};
\begin{scope}[xshift=0cm]
\begin{scope}[yshift=0cm]
\draw[retained] (10.2,0) arc (-90:90:0.6) -- (0,1.2) arc (90:270:0.6);
\draw[retained] (0,0) -- (0.55,0);
\coordinate (before-0-f1-u) at (0.55,0);
\coordinate (before-0-f1-v) at (1.65,0);
\draw[removed] (0.55,0) -- node[midway,below=2pt] {$\phi(f_{1})$} (1.65,0);
\draw[retained] (1.65,0) -- (10.2,0);
\node at (5.1,0.6) {$C(0)$};
\end{scope}
\begin{scope}[yshift=-2.8cm]
\draw[retained] (10.2,0) arc (-90:90:0.6) -- (0,1.2) arc (90:270:0.6);
\draw[retained] (0,0) -- (0.55,0);
\coordinate (before-1-f1-u) at (0.55,0);
\coordinate (before-1-f1-v) at (1.65,0);
\draw[removed] (0.55,0) -- node[midway,below=2pt] {$\phi(f_{1})$} (1.65,0);
\draw[retained] (1.65,0) -- (3.1,0);
\coordinate (before-1-f2-u) at (3.1,0);
\coordinate (before-1-f2-v) at (4.2,0);
\draw[removed] (3.1,0) -- node[midway,below=2pt] {$\phi(f_{2})$} (4.2,0);
\draw[retained] (4.2,0) -- (10.2,0);
\node at (5.1,0.6) {$C(1)$};
\end{scope}
\begin{scope}[yshift=-6.8cm]
\draw[retained] (10.2,0) arc (-90:90:0.6) -- (0,1.2) arc (90:270:0.6);
\draw[retained] (0,0) -- (5.8,0);
\coordinate (before-2-fp-2-u) at (5.8,0);
\coordinate (before-2-fp-2-v) at (6.9,0);
\draw[removed] (5.8,0) -- node[midway,below=2pt] {$\phi(f_{p-2})$} (6.9,0);
\draw[retained] (6.9,0) -- (8.35,0);
\coordinate (before-2-fp-1-u) at (8.35,0);
\coordinate (before-2-fp-1-v) at (9.45,0);
\draw[removed] (8.35,0) -- node[midway,below=2pt] {$\phi(f_{p-1})$} (9.45,0);
\draw[retained] (9.45,0) -- (10.2,0);
\node at (5.1,0.6) {$C(p-2)$};
\end{scope}
\begin{scope}[yshift=-9.6cm]
\draw[retained] (10.2,0) arc (-90:90:0.6) -- (0,1.2) arc (90:270:0.6);
\draw[retained] (0,0) -- (8.35,0);
\coordinate (before-3-fp-1-u) at (8.35,0);
\coordinate (before-3-fp-1-v) at (9.45,0);
\draw[removed] (8.35,0) -- node[midway,below=2pt] {$\phi(f_{p-1})$} (9.45,0);
\draw[retained] (9.45,0) -- (10.2,0);
\node at (5.1,0.6) {$C(p-1)$};
\end{scope}
\node[font=\Large] at (5.1,-4.5) {$\vdots$};
\end{scope}
\begin{scope}[xshift=14cm]
\begin{scope}[yshift=0cm]
\draw[retained] (10.2,0) arc (-90:90:0.6) -- (0,1.2) arc (90:270:0.6);
\draw[retained] (0,0) -- (0.55,0);
\coordinate (after-0-f1-u) at (0.55,0);
\coordinate (after-0-f1-v) at (1.65,0);
\draw[retained] (1.65,0) -- (10.2,0);
\end{scope}
\begin{scope}[yshift=-2.8cm]
\draw[retained] (10.2,0) arc (-90:90:0.6) -- (0,1.2) arc (90:270:0.6);
\draw[retained] (0,0) -- (0.55,0);
\coordinate (after-1-f1-u) at (0.55,0);
\coordinate (after-1-f1-v) at (1.65,0);
\draw[retained] (1.65,0) -- (3.1,0);
\coordinate (after-1-f2-u) at (3.1,0);
\coordinate (after-1-f2-v) at (4.2,0);
\draw[retained] (4.2,0) -- (10.2,0);
\end{scope}
\begin{scope}[yshift=-6.8cm]
\draw[retained] (10.2,0) arc (-90:90:0.6) -- (0,1.2) arc (90:270:0.6);
\draw[retained] (0,0) -- (5.8,0);
\coordinate (after-2-fp-2-u) at (5.8,0);
\coordinate (after-2-fp-2-v) at (6.9,0);
\draw[retained] (6.9,0) -- (8.35,0);
\coordinate (after-2-fp-1-u) at (8.35,0);
\coordinate (after-2-fp-1-v) at (9.45,0);
\draw[retained] (9.45,0) -- (10.2,0);
\end{scope}
\begin{scope}[yshift=-9.6cm]
\draw[retained] (10.2,0) arc (-90:90:0.6) -- (0,1.2) arc (90:270:0.6);
\draw[retained] (0,0) -- (8.35,0);
\coordinate (after-3-fp-1-u) at (8.35,0);
\coordinate (after-3-fp-1-v) at (9.45,0);
\draw[retained] (9.45,0) -- (10.2,0);
\end{scope}
\node[font=\Large] at (5.1,-4.5) {$\vdots$};
\end{scope}
\draw[->,line width=1.1pt] (11.65,-4.2) -- (12.85,-4.2);
% First joining pair: the copies of f_1 in C(0) and C(1).
\draw[joined] (after-0-f1-u) -- (after-1-f1-u);
\draw[joined] (after-0-f1-v) -- (after-1-f1-v);
% The next pair continues into the omitted rows.
\draw[joined] (after-1-f2-u) -- (17.1,-4.25);
\draw[joined] (after-1-f2-v) -- (18.2,-4.25);
% The omitted rows continue into the copies of f_{p-2} in C(p-2).
\draw[joined] (19.8,-5.0) -- (after-2-fp-2-u);
\draw[joined] (20.9,-5.0) -- (after-2-fp-2-v);
% Last joining pair: the copies of f_{p-1} in C(p-2) and C(p-1).
\draw[joined] (after-2-fp-1-u) -- (after-3-fp-1-u);
\draw[joined] (after-2-fp-1-v) -- (after-3-fp-1-v);
\node[dot] at (before-0-f1-u) {};
\node[dot] at (before-0-f1-v) {};
\node[dot] at (before-1-f1-u) {};
\node[dot] at (before-1-f1-v) {};
\node[dot] at (before-1-f2-u) {};
\node[dot] at (before-1-f2-v) {};
\node[dot] at (before-2-fp-2-u) {};
\node[dot] at (before-2-fp-2-v) {};
\node[dot] at (before-2-fp-1-u) {};
\node[dot] at (before-2-fp-1-v) {};
\node[dot] at (before-3-fp-1-u) {};
\node[dot] at (before-3-fp-1-v) {};
\node[dot] at (after-0-f1-u) {};
\node[dot] at (after-0-f1-v) {};
\node[dot] at (after-1-f1-u) {};
\node[dot] at (after-1-f1-v) {};
\node[dot] at (after-1-f2-u) {};
\node[dot] at (after-1-f2-v) {};
\node[dot] at (after-2-fp-2-u) {};
\node[dot] at (after-2-fp-2-v) {};
\node[dot] at (after-2-fp-1-u) {};
\node[dot] at (after-2-fp-1-v) {};
\node[dot] at (after-3-fp-1-u) {};
\node[dot] at (after-3-fp-1-v) {};
\draw[removed] (2.5,-11.5) -- (3.35,-11.5);
\node[anchor=west] at (3.5,-11.5) {deleted edges from $\phi(f_j)$};
\draw[joined] (12.5,-11.5) -- (13.35,-11.5);
\node[anchor=west] at (13.5,-11.5) {added edges, each of difference $1$};
\end{tikzpicture}
\endgroup
\caption{Joining the cycles $C(0),C(1),\ldots,C(p-1)$ in Type~2 of the proof of Theorem~\ref{thm:q}.}
\label{fig:type2}
\end{center}
\vspace{-0.5cm}
\end{figure}

\smallskip\noindent\textbf{Type~3:} $\alpha(C^{\ell})=0$ and $n\geq 3$.\smallskip

By Lemma~\ref{le:cycle_lift}, $(C^{\ell})^{\alpha}$ consists of $q$ vertex-disjoint cycles of length $2q$.
We shall modify these cycles into a Hamilton cycle of $H_q$ while avoiding the fibres already used for Type~1.
Since $q=4\ell-1$ in this case, we have $\frac{q-1}{2}=2\ell -1$.
Let $K=\{1, 2, \ldots, 2\ell-1\}\setminus \{\ell\}$.\footnote{In this paragraph, sets of vertex indices and all operations on them are interpreted in $\mathbb{Z}/q\mathbb{Z}$, with the residue $0$ represented by $q$.}
For $k\neq \ell$, the modification for Type~1 uses only fibres above the end-vertices of $e_k=u_{q-k}v_k=u_{-k}v_{k}$.
Thus the vertices in $U$ whose fibres are already used have indices in $-K$, and those in $V$ have indices in $K$.
Consider the perfect matching $M_{2\ell+1}=\{u_iv_{i+2\ell}\colon\, 1\leq i\leq q\}$.
The edge $u_iv_{i+2\ell}$ has an end-vertex whose fibre was used in Type~1 precisely when $i\in -K$ or $i+2\ell\in K$.
Since replacing $j\in K$ by $2\ell-j$ permutes $K$, we have $2\ell-K=K$, i.e., $K-2\ell=-K$, so $i+2\ell\in K$ is equivalent to $i\in -K$.
Hence, exactly $|K|=2\ell-2$ edges of $M_{2\ell+1}$ meet a previously used fibre.
The number of remaining edges in $M_{2\ell+1}$ is $q-(2\ell-2)=\frac{q+3}{2}$.
Since $q=p^n$ and $n\geq 3$, we have $\frac{q+3}{2}\geq p$.
Choose $p$ edges in $M_{2\ell+1}$ such that the fibres of their end-vertices were not used in Type~1, and denote them by $e, f'_1, \ldots, f'_{p-1}$.

Since $p$ is odd, we have $2\neq 0$.
By Lemma~\ref{le:difference}, we have $S_U\cup S_V=\mathbb{F}_q^{\times}$, so $2$ belongs to at least one of $S_U$ and $S_V$.
Denote the end-vertices of $e$ by $w$ and $t$, choosing $w$ so that the difference set of $\phi(w)$ contains $2$.
Let $\beta$ be the voltage of the arc from $w$ to $t$.
Orient $C^{\ell}$ so that it traverses $e$ from $w$ to $t$.
For each $x\in \mathbb{F}_q$, let $C(x)$ be the component of $(C^{\ell})^\alpha$ containing the arc $(e, x)$ from $(w,x)$ to $(t, x+\beta)$.
Let $L=\mathbb{F}_p$.
Recall that the cosets of $L$ form a partition of $\mathbb{F}_q$.
We now modify the cycles $C(x)$ ($x\in \mathbb{F}_q$) into a Hamilton cycle of $H_q$ in two steps.

\smallskip\noindent\textit{Step~1. Join the components with indices within each coset of $L$.} \smallskip

By Lemma~\ref{le:difference}, the difference set of $\phi(t)$ is symmetric and its images span $\mathbb{F}_q/L$.
Apply Lemma~\ref{le:z_i} with $a=1$ to obtain a cyclic ordering $A_0,A_1,\ldots,A_{m-1},A_m=A_0$ of the cosets of $L$
and elements $z_0,z_1, \ldots,z_{m-1}$ in the difference set of $\phi(t)$ satisfying~\eqref{eq:splice}, where $m=\frac{q}{p}=p^{n-1}$.

Take an $x_0\in A_0$ and define successively $x_{i+1}\colonequals x_i+z_i-2$ for $0\leq i\leq m-1$.
By Observation~\ref{obs:coset}, we have $x_i\in A_i$ for $0\leq i\leq m-1$, and $x_m=x_0$.
Then for each $0\leq i\leq m-1$, we have $A_i=\{x_i, x_i+1, \ldots, x_i+p-1\}$.
Let $\mathcal{D}_i\colonequals \{C(x)\colon\, x\in A_i\}=\{C(x_i), C(x_i+1), \ldots, C(x_i+p-1)\}$.
Then $\mathcal{D}_0, \mathcal{D}_1, \ldots, \mathcal{D}_{m-1}$ partition the $q$ components of $(C^{\ell})^{\alpha}$ into $m$ parts each containing $p$ components.

Fix an $i\in \{0, 1, \ldots, m-1\}$.
We join the cycles $C(x_i), C(x_i+1), \ldots, C(x_i+p-1)$ in the same way as in Type~2.
Start with $C(x_i)$.
In step $j$, delete the lifted copy of $f'_j$ in the original component $C(x_i+j-1)$ from the cycle already constructed, and delete the lifted copy of $f'_j$ in $C(x_i+j)$.
Both copies are still present since the base edges used in earlier steps are distinct from $f'_j$.
We obtain two paths, and then join their end-vertices using two edges inside the fibres above the end-vertices of $f'_j$.
Such edges exist since $1\in S_U\cap S_V$ by Lemma~\ref{le:difference}.
Then we obtain a single cycle on $V\big(C(x_i)\big)\cup \cdots \cup V\big(C(x_i+j)\big)$.
After $p-1$ steps, we obtain a cycle $D(i)$ containing all vertices of the $p$ components in $\mathcal{D}_i$.
For different $j$, the added edges lie in distinct vertex fibres since $f'_1, f'_2, \ldots, f'_{p-1}$ belong to a matching.
For a fixed $j$, different values of $i$ use disjoint sets of original components, so their added edges have disjoint end-vertices.
Hence, all added edges are distinct.

\smallskip\noindent\textit{Step~2. Join the cycles $D(0), D(1), \ldots, D(m-1)$.} \smallskip

We orient each $C(x)$ by the orientation inherited from $C^\ell$.
At step~$j$ of Step~1, we join corresponding end-vertices of the two deleted lifted copies of $f'_j$.
Hence, in the resulting cycle, the surviving edges of $C(x_i+j)$ are traversed in the opposite direction to those of $C(x_i+j-1)$, relative to their inherited orientations.
The lifted copies of $e$ survive every joining step since $e$ is distinct from every $f'_j$.
It follows inductively that the copies $(e,x_i)$ and $(e,x_i+2)$ are traversed in the same direction in $D(i)$.

We join the cycles $D(0), D(1), \ldots, D(m-1)$ as in the proof of Lemma~\ref{le:join}.
For each $i\in \{0,1,\ldots, m-1\}$, we delete $(e,x_i)$ and $(e,x_i+2)$ from $D(i)$.
Then we obtain two paths, one with end-vertices $(t, x_i+\beta)$ and $(w, x_i+2)$, and the other one with end-vertices $(t, x_i+2+\beta)$ and $(w, x_i)$.
Now we add the edge between $(w, x_i)$ and $(w, x_i+2)$, which exists since the difference set of $\phi(w)$ contains $2$.
Since  $(e,x_i)$ and $(e,x_i+2)$ have the same direction, this joins the two paths into one path $P(i)$ with end-vertices $(t, x_i+\beta)$ and $(t, x_i+2+\beta)$, and $P(i)$ contains every vertex of $D(i)$.
It remains to join the $m$ paths $P(0), P(1), \ldots, P(m-1)$ into a Hamilton cycle of $H_q$.
For every $i$ (with indices modulo $m$), add an edge from $(t,x_i+\beta)$ to $(t,x_{i+1}+2+\beta)$.
Such an edge exists since the difference between its second coordinates is $x_{i+1}+2-x_i=z_i$ and $z_i$ belongs to the difference set of $\phi(t)$.
This joins the paths in the cyclic order $P(0), P(1), \ldots, P(m-1), P(0)$.
The resulting graph is a cycle containing all vertices of all the $D(i)$ (and thus every vertex of $H_q$), so it is a Hamilton cycle of $H_q$.

We now verify edge-disjointness.
Step~1 adds edges only inside the fibres above the end-vertices of $f'_1, \ldots, f'_{p-1}$, and Step~2 adds edges only inside the fibres above the end-vertices of $e$.
Since $e, f'_1, \ldots, f'_{p-1}$ belong to one perfect matching, the fibres above their end-vertices are distinct.
Moreover, by the choice of $e, f'_1, \ldots, f'_{p-1}$, none of these fibres was used for Type~1.
Furthermore, since $C^1, \ldots, C^{(q-1)/2}$ are edge-disjoint in the base graph $G_q$, all the cross edges of $(C^1)^{\alpha}, \ldots, (C^{(q-1)/2})^{\alpha}$ are distinct.
Hence, the Hamilton cycle just constructed is edge-disjoint from every Hamilton cycle obtained from Type~1.

Combining Types~1, 2 and 3, we obtain $\frac{q-1}{2}$ edge-disjoint Hamilton cycles of $H_q$ for every odd prime power $q\geq 3$.
This completes the proof for odd $q$, and thus completes the proof of Theorem~\ref{thm:q}.
\end{proof}

\section*{Declaration of competing interests}
\vspace{-0.2cm}

The author declares that they have no competing interests.
\vspace{-0.2cm}

\section*{Data availability}
\vspace{-0.2cm}

No data were used for the research described in the article.
\vspace{-0.2cm}

\begin{spacing}{0.8} %行间距

\end{spacing}

%\newpage
\appendix
\section*{Appendix}

\section{Constructing the Hamilton cycles using $O(|E(H_q)|)$ operations}
\label{sec:app-2}

In this appendix, we show that the $\big\lceil\frac{q-1}{2}\big\rceil$ edge-disjoint Hamilton cycles given in Theorem~\ref{thm:q} can be constructed using $O(|E(H_q)|)$ finite field operations, comparisons, and elementary list operations.\footnote{Proposition~\ref{prop:time} is an operation-count statement. It does not assert an $O(|E(H_q)|)$ bound on bit operations. Indeed, there are $\lceil(q-1)/2\rceil$ output cycles, each having $2q^2$ edges. In the usual explicit binary representation, a vertex label requires $\Theta(\log q)$ bits. Thus the output alone has $\Theta(q^3\log q)$ bits, whereas $|E(H_q)|=\Theta(q^3)$.}

\begin{proposition}\label{prop:time}
Let $q\geq3$ be a prime power.
Given the finite field $\mathbb{F}_q$ and the primitive element $\xi$ used in the definition of $H_q$,
the construction in Theorem~\ref{thm:q} produces $\big\lceil\frac{q-1}{2}\big\rceil$ edge-disjoint Hamilton cycles
using $O(|E(H_q)|)$ finite field operations, comparisons, and elementary list operations.
\end{proposition}

\begin{proof}
Write $q=p^n=4\ell+\delta$, where $\delta\in\{-1,0,1\}$.
Recall that $H_q$ has $2q^2$ vertices and is $\frac{3q-\delta}{2}$-regular, so $|E(H_q)|=q^2\cdot\frac{3q-\delta}{2}=\Theta(q^3)$.
We shall show that each of the $\big\lceil\frac{q-1}{2}\big\rceil=O(q)$ Hamilton cycles can be constructed using $O(q^2)$ operations.
Each Hamilton cycle is output as a list of its edges.
The operations counted below are addition, subtraction, and multiplication in $\mathbb{F}_q$, comparison of labels, and elementary list operations such as reading the next entry and appending an entry.

We first generate the powers of the primitive element and the difference sets $S_U$ and $S_V$.
Using the explicit descriptions in the proof of Lemma~\ref{le:difference}, these lists are generated using $O(q)$ field operations and list operations.
The remaining construction has three parts.

\smallskip\noindent
\emph{1. Constructing the cyclic ordering of cosets and the elements $z_i$ in Lemma~\ref{le:z_i}.}
Let $L=a\mathbb{F}_p$ with $a\neq 0$, and let $S=-S\subseteq \mathbb{F}_q^{\times}$ have images spanning $\mathbb{F}_q/L$.
We can generate the initial list $L$ by repeated addition of $a$, using $O(p)$ operations.

First, we choose $b_1,\ldots,b_{n-1}\in S$ such that the cosets $b_1+L, \ldots, b_{n-1}+L$ form a basis of $\mathbb{F}_q/L$.
Start with $W=L$ and scan the elements of $S$.
Whenever $b\notin W$, retain $b$ and replace $W$ by $W+\mathbb{F}_pb$.
The condition $b\notin W$ means that $b+L$ is not in the span of the previously retained cosets.
Thus the retained cosets are linearly independent.
Since the images of $S$ span $\mathbb{F}_q/L$, the procedure reaches $W=\mathbb{F}_q$ after retaining exactly $n-1$ elements, whose images form the required basis.
For the operation count, we store $W$ as a list.
Each membership test uses at most $q$ comparisons, and there are at most $q$ tests.
When retaining $b$, construct the enlarged list by appending $W+b,W+2b,\ldots,W+(p-1)b$ to the old $W$.
These cosets are disjoint, so this requires at most $O(q)$ operations per enlargement.
There are $n-1\leq q$ enlargements.
Hence, the entire procedure uses $O(q^2)$ operations to choose the required $b_1,\ldots,b_{n-1}\in S$.

If $q=4$, the selected element $b_1\in S\setminus L$ gives the ordering $L,b_1+L,L$ and the differences $z_0=z_1=b_1$ directly, using a fixed number of operations.
For $q\neq 4$, we proceed as follows.
The construction of the Hamilton cycle in Lemma~\ref{le:Cayley} uses $O(p^d)$ operations, where $d$ is the dimension of the vector space.
Indeed, the initial cycle requires $O(p)$ operations for odd $p$, and $O(1)$ operations for $p=2$.
For the inductive step, at each subsequent stage $j$, we form $p$ paths each with $p^{j-1}$ vertices and concatenate their vertex lists.
This takes $O(p^j)$ operations.
Summing over all stages, this takes $O(p^d)$ operations.
For Lemma~\ref{le:z_i}, we apply this with $d=n-1$, so we can construct the cyclic vertex list in $O(m)$ operations, where $m=p^{n-1}=\frac{q}{p}$.
Computing $z_i$ for $0\leq i\leq m-1$ requires another $O(m)$ operations.

Therefore, the cyclic ordering of cosets and the elements $z_i$ in Lemma~\ref{le:z_i} can be constructed in $O(q^2)+O(m)=O(q^2)$ operations.

\smallskip\noindent
\emph{2. Listing the required edge modifications.}
For each base cycle $C^k$, its cyclic vertex order in~\eqref{eq:Ck}, its arc voltages, and its net voltage $\alpha(C^k)$ can be generated using $O(q)$ operations.
We then construct two lists: the cross edges to be deleted from $(C^k)^\alpha$, and the edges within fibres of vertices to be added.

First, suppose that $\alpha(C^k)\neq 0$.
This corresponds to the case when $q$ is even and Type~1 in the case when $q$ is odd.
If $n=1$, no modification is needed.
For $n\geq 2$, we use Corollary~\ref{cor:join} (in fact, Lemma~\ref{le:join}) to make four edge changes for each of the $\frac{q}{p}$ cosets.
The required conditions are guaranteed by Lemma~\ref{le:z_i}, and we use $O(q^2)$ operations for Lemma~\ref{le:z_i} as shown above.
Moreover, for every $k$, the deletion and addition lists have length $O(q)$.

Next, suppose that $\alpha(C^k)= 0$.
This corresponds to Types~2 and~3 in the case when $q$ is odd.
In Type~2, we join the cycles in $p-1$ steps.
In Type~3, choosing $e,f'_1,\ldots,f'_{p-1}$ takes $O(q^2)$ operations by scanning $M_{2\ell+1}$ and checking whether its edges meet fibres already used for Type~1.
In Step~1 of Type~3, we join the cycles in $\frac{q}{p}(p-1)<q$ steps, and Step~2 makes $O(\frac{q}{p})$ edge changes.
Each required edge in an original component $C(x)$ can be located by following the lift of $C^{\ell}$ from the vertex used to index that component.
Since $C(x)$ has length $2q$, this takes $O(q)$ operations per edge.
In Type 3, we also use Lemma~\ref{le:z_i}, which costs $O(q^2)$ operations.
Thus all the required edges can be listed using $O(q^2)$ operations.

Therefore, for each $C^k$, there are $O(q)$ deleted and added edges, and their lists can be constructed using $O(q^2)$ operations.

\smallskip\noindent
\emph{3. Outputting the Hamilton cycles.}
For each $C^k$, generate the lifted edges in order by traversing the $2q$ base edges in~\eqref{eq:Ck}.
For each base edge, oriented as in~\eqref{eq:Ck}, enumerate its $q$ lifted edges in a fixed order of their field coordinates (i.e., the second coordinates).
This uses $O(q^2)$ operations.
For each deleted edge, scan the list~\eqref{eq:Ck} to find the position and orientation of its base edge.
Order its two end-vertices according to that orientation.
Scan the fixed list of field elements to find the position of its initial field coordinate.
There are $O(q)$ deleted edges, and each of these scans uses $O(q)$ operations, so this requires $O(q^2)$ operations in total.
Sort the deletion list lexicographically by the recorded base edge position and field coordinate position, using merge sort with $O(q\log q)$ comparisons and list operations.
Scanning the two lists together then outputs all undeleted edges in $O(q^2)$ operations.
Finally, append the $O(q)$ added edges within fibres of vertices.
Thus outputting one Hamilton cycle requires $O(q^2)$ operations.

\smallskip
Combining the above three parts, each Hamilton cycle requires $O(q^2)$ operations.
Since there are $\big\lceil\frac{q-1}{2}\big\rceil=O(q)$ Hamilton cycles, the total is $O(q^3)=O(|E(H_q)|)$ operations.
\end{proof}

%\vspace{-0.5cm}

\section{Illustrations of $G_q$ and $H_q$ for $q\in \{3,4\}$}
\label{sec:app-1}

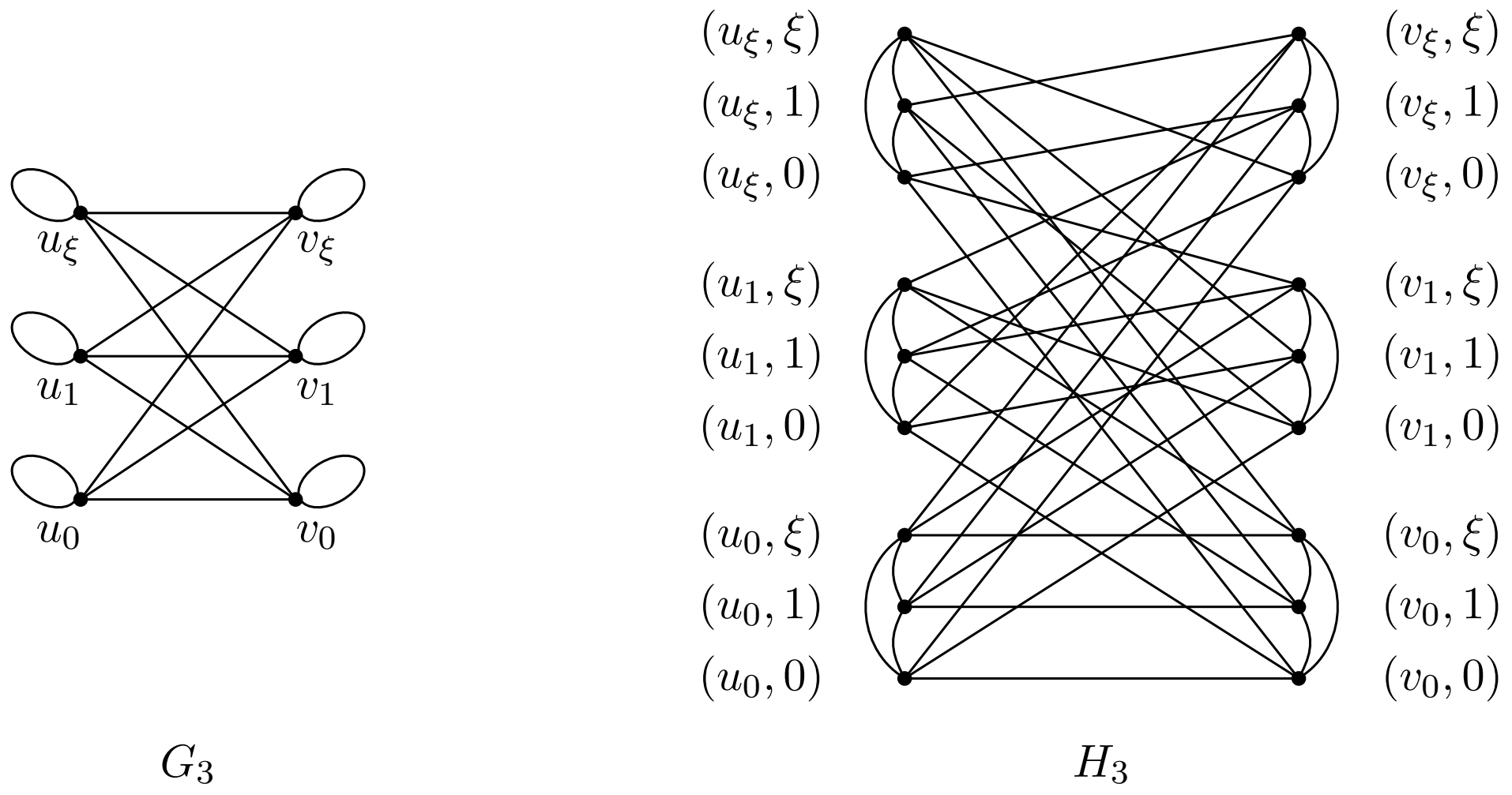
\begin{figure}[htbp]
\begin{center}
\begin{tikzpicture}[scale=0.06,auto,swap]
\tikzstyle{vertex1}=[circle,draw=black,fill=black]

\node[vertex1,scale=0.3] (x10) at (10,40) {}; \draw (7,35) node {$u_{\xi}$}; \draw[line width=0.2mm, rotate around={-30:(5,42.5)}] (5,42.5) ellipse (5 and 3);
\node[vertex1,scale=0.3] (x20) at (10,20) {}; \draw (7,15) node {$u_1$}; \draw[line width=0.2mm, rotate around={-30:(5,22.5)}] (5,22.5) ellipse (5 and 3);
\node[vertex1,scale=0.3] (x30) at (10,0) {}; \draw (7,-5) node {$u_0$}; \draw[line width=0.2mm, rotate around={-30:(5,2.5)}] (5,2.5) ellipse (5 and 3);
\node[vertex1,scale=0.3] (x11) at (40,40) {}; \draw (43,35) node {$v_{\xi}$}; \draw[line width=0.2mm, rotate around={30:(45,42.5)}] (45,42.5) ellipse (5 and 3);
\node[vertex1,scale=0.3] (x21) at (40,20) {}; \draw (43,15) node {$v_1$}; \draw[line width=0.2mm, rotate around={30:(45,22.5)}] (45,22.5) ellipse (5 and 3);
\node[vertex1,scale=0.3] (x31) at (40,0) {}; \draw (43,-5) node {$v_0$}; \draw[line width=0.2mm, rotate around={30:(45,2.5)}] (45,2.5) ellipse (5 and 3);
\draw[line width=0.2mm, black, opacity=1] (x10) -- (x11); \draw[line width=0.2mm, black, opacity=1] (x10) -- (x21); \draw[line width=0.2mm, black, opacity=1] (x10) -- (x31);
\draw[line width=0.2mm, black, opacity=1] (x20) -- (x11); \draw[line width=0.2mm, black, opacity=1] (x20) -- (x21); \draw[line width=0.2mm, black, opacity=1] (x20) -- (x31);
\draw[line width=0.2mm, black, opacity=1] (x30) -- (x11); \draw[line width=0.2mm, black, opacity=1] (x30) -- (x21); \draw[line width=0.2mm, black, opacity=1] (x30) -- (x31);
\draw (25,-37) node {$G_3$};

\node[vertex1,scale=0.3] (y101) at (125,65) {}; \draw (105,65) node {$(u_{\xi},\xi)$};
\node[vertex1,scale=0.3] (y102) at (125,55) {}; \draw (105,55) node {$(u_{\xi},1)$};
\node[vertex1,scale=0.3] (y103) at (125,45) {}; \draw (105,45) node {$(u_{\xi},0)$};
\draw[line width=0.2mm, black, opacity=1] (y101) .. controls (123,61) and (123,59) .. (y102);
\draw[line width=0.2mm, black, opacity=1] (y102) .. controls (123,51) and (123,49) .. (y103);
\draw[line width=0.2mm, black, opacity=1] (y101) .. controls (118,60) and (118,50) .. (y103);
\node[vertex1,scale=0.3] (y201) at (125,30) {}; \draw (105,30) node {$(u_1,\xi)$};
\node[vertex1,scale=0.3] (y202) at (125,20) {}; \draw (105,20) node {$(u_1,1)$};
\node[vertex1,scale=0.3] (y203) at (125,10) {}; \draw (105,10) node {$(u_1,0)$};
\draw[line width=0.2mm, black, opacity=1] (y201) .. controls (123,26) and (123,24) .. (y202);
\draw[line width=0.2mm, black, opacity=1] (y202) .. controls (123,16) and (123,14) .. (y203);
\draw[line width=0.2mm, black, opacity=1] (y201) .. controls (118,25) and (118,15) .. (y203);
\node[vertex1,scale=0.3] (y301) at (125,-5) {}; \draw (105,-5) node {$(u_0,\xi)$};
\node[vertex1,scale=0.3] (y302) at (125,-15) {}; \draw (105,-15) node {$(u_0,1)$};
\node[vertex1,scale=0.3] (y303) at (125,-25) {}; \draw (105,-25) node {$(u_0,0)$};
\draw[line width=0.2mm, black, opacity=1] (y301) .. controls (123,-9) and (123,-11) .. (y302);
\draw[line width=0.2mm, black, opacity=1] (y302) .. controls (123,-19) and (123,-21) .. (y303);
\draw[line width=0.2mm, black, opacity=1] (y301) .. controls (118,-10) and (118,-20) .. (y303);
\node[vertex1,scale=0.3] (y111) at (180,65) {}; \draw (200,65) node {$(v_{\xi},\xi)$};
\node[vertex1,scale=0.3] (y112) at (180,55) {}; \draw (200,55) node {$(v_{\xi},1)$};
\node[vertex1,scale=0.3] (y113) at (180,45) {}; \draw (200,45) node {$(v_{\xi},0)$};
\draw[line width=0.2mm, black, opacity=1] (y111) .. controls (182,61) and (182,59) .. (y112);
\draw[line width=0.2mm, black, opacity=1] (y112) .. controls (182,51) and (182,49) .. (y113);
\draw[line width=0.2mm, black, opacity=1] (y111) .. controls (187,60) and (187,50) .. (y113);
\node[vertex1,scale=0.3] (y211) at (180,30) {}; \draw (200,30) node {$(v_1,\xi)$};
\node[vertex1,scale=0.3] (y212) at (180,20) {}; \draw (200,20) node {$(v_1,1)$};
\node[vertex1,scale=0.3] (y213) at (180,10) {}; \draw (200,10) node {$(v_1,0)$};
\draw[line width=0.2mm, black, opacity=1] (y211) .. controls (182,26) and (182,24) .. (y212);
\draw[line width=0.2mm, black, opacity=1] (y212) .. controls (182,16) and (182,14) .. (y213);
\draw[line width=0.2mm, black, opacity=1] (y211) .. controls (187,25) and (187,15) .. (y213);
\node[vertex1,scale=0.3] (y311) at (180,-5) {}; \draw (200,-5) node {$(v_0,\xi)$};
\node[vertex1,scale=0.3] (y312) at (180,-15) {}; \draw (200,-15) node {$(v_0,1)$};
\node[vertex1,scale=0.3] (y313) at (180,-25) {}; \draw (200,-25) node {$(v_0,0)$};
\draw[line width=0.2mm, black, opacity=1] (y311) .. controls (182,-9) and (182,-11) .. (y312);
\draw[line width=0.2mm, black, opacity=1] (y312) .. controls (182,-19) and (182,-21) .. (y313);
\draw[line width=0.2mm, black, opacity=1] (y311) .. controls (187,-10) and (187,-20) .. (y313);
\draw[line width=0.2mm, black, opacity=1] (y101) -- (y113); \draw[line width=0.2mm, black, opacity=1] (y102) -- (y111); \draw[line width=0.2mm, black, opacity=1] (y103) -- (y112);
\draw[line width=0.2mm, black, opacity=1] (y101) -- (y212); \draw[line width=0.2mm, black, opacity=1] (y102) -- (y213); \draw[line width=0.2mm, black, opacity=1] (y103) -- (y211);
\draw[line width=0.2mm, black, opacity=1] (y101) -- (y311); \draw[line width=0.2mm, black, opacity=1] (y102) -- (y312); \draw[line width=0.2mm, black, opacity=1] (y103) -- (y313);
\draw[line width=0.2mm, black, opacity=1] (y201) -- (y112); \draw[line width=0.2mm, black, opacity=1] (y202) -- (y113); \draw[line width=0.2mm, black, opacity=1] (y203) -- (y111);
\draw[line width=0.2mm, black, opacity=1] (y201) -- (y213); \draw[line width=0.2mm, black, opacity=1] (y202) -- (y211); \draw[line width=0.2mm, black, opacity=1] (y203) -- (y212);
\draw[line width=0.2mm, black, opacity=1] (y201) -- (y311); \draw[line width=0.2mm, black, opacity=1] (y202) -- (y312); \draw[line width=0.2mm, black, opacity=1] (y203) -- (y313);
\draw[line width=0.2mm, black, opacity=1] (y301) -- (y111); \draw[line width=0.2mm, black, opacity=1] (y302) -- (y112); \draw[line width=0.2mm, black, opacity=1] (y303) -- (y113);
\draw[line width=0.2mm, black, opacity=1] (y301) -- (y211); \draw[line width=0.2mm, black, opacity=1] (y302) -- (y212); \draw[line width=0.2mm, black, opacity=1] (y303) -- (y213);
\draw[line width=0.2mm, black, opacity=1] (y301) -- (y311); \draw[line width=0.2mm, black, opacity=1] (y302) -- (y312); \draw[line width=0.2mm, black, opacity=1] (y303) -- (y313);
\draw (152.5,-37) node {$H_3$};		
\end{tikzpicture}
\caption{The base graph $G_3$ and the MMS graph $H_3$. In this case, $\ell=1$, $\delta=-1$, $\mathbb{F}_3=\{\xi(=2), 1, 0\}$ and $Y=Y'=\{1\}$.}
\label{fig:q=3}
\end{center}
\vspace{-0.5cm}
\end{figure}

We illustrate the base graphs $G_q$ and their lifts $H_q$ for $q=3$ and $q=4$ in Figures~\ref{fig:q=3} and~\ref{fig:q=4}, respectively.
In this appendix, the subscripts of $u_x$ and $v_y$ are field elements, as in Section~\ref{subsec:H_q}, rather than the integer indices introduced in Section~\ref{subsec:cycles}.
%\vspace{0.5cm}

\begin{figure}[htbp]
\begin{center}
\definecolor{MMSblue}{RGB}{0,90,165}
\definecolor{MMSorange}{RGB}{190,74,0}
\definecolor{MMSgreen}{RGB}{0,125,100}
\definecolor{MMSpurple}{RGB}{135,65,150}
\begin{tikzpicture}[scale=0.6,x=1cm,y=1cm,font=\small]
\tikzset{dot/.style={circle,fill=black,draw=black,inner sep=0pt,minimum size=3pt}}
% Move only the base graph to the left of the first row.
\begin{scope}[shift={(-9,0)}]
\node at (11,14.5) {(a) $G_4$};
% G_4 is K_{4,4} with two semiedges at each vertex.
\coordinate (g02) at (8,18.5);
\coordinate (g12) at (14,18.5);
\coordinate (g03) at (8,17.5);
\coordinate (g13) at (14,17.5);
\coordinate (g01) at (8,16.5);
\coordinate (g11) at (14,16.5);
\coordinate (g00) at (8,15.5);
\coordinate (g10) at (14,15.5);
% The four edges of M_1; their lifts appear in the same colour below.
\draw[draw=MMSblue,line width=0.65pt] (g02) -- (g12);
\draw[draw=MMSblue,line width=0.65pt] (g03) -- (g13);
\draw[draw=MMSblue,line width=0.65pt] (g01) -- (g11);
\draw[draw=MMSblue,line width=0.65pt] (g00) -- (g10);
% The four edges of M_2; their lifts appear in the same colour below.
\draw[draw=MMSorange,line width=0.65pt] (g02) -- (g13);
\draw[draw=MMSorange,line width=0.65pt] (g03) -- (g11);
\draw[draw=MMSorange,line width=0.65pt] (g01) -- (g10);
\draw[draw=MMSorange,line width=0.65pt] (g00) -- (g12);
% The four edges of M_3; their lifts appear in the same colour below.
\draw[draw=MMSgreen,line width=0.65pt] (g02) -- (g11);
\draw[draw=MMSgreen,line width=0.65pt] (g03) -- (g10);
\draw[draw=MMSgreen,line width=0.65pt] (g01) -- (g12);
\draw[draw=MMSgreen,line width=0.65pt] (g00) -- (g13);
% The four edges of M_4; their lifts appear in the same colour below.
\draw[draw=MMSpurple,line width=0.65pt] (g02) -- (g10);
\draw[draw=MMSpurple,line width=0.65pt] (g03) -- (g12);
\draw[draw=MMSpurple,line width=0.65pt] (g01) -- (g13);
\draw[draw=MMSpurple,line width=0.65pt] (g00) -- (g11);
\draw[black,line width=0.55pt] (g02) -- (7.45,18.68);
\draw[black,line width=0.55pt] (g02) -- (7.45,18.32);
\node[anchor=east,inner sep=0pt] at (7.15,18.5) {$u_{\xi}$};
\node[dot] at (g02) {};
\draw[black,line width=0.55pt] (g12) -- (14.55,18.68);
\draw[black,line width=0.55pt] (g12) -- (14.55,18.32);
\node[anchor=west,inner sep=0pt] at (14.85,18.5) {$v_{\xi}$};
\node[dot] at (g12) {};
\draw[black,line width=0.55pt] (g03) -- (7.45,17.68);
\draw[black,line width=0.55pt] (g03) -- (7.45,17.32);
\node[anchor=east,inner sep=0pt] at (7.15,17.5) {$u_{\xi^2}$};
\node[dot] at (g03) {};
\draw[black,line width=0.55pt] (g13) -- (14.55,17.68);
\draw[black,line width=0.55pt] (g13) -- (14.55,17.32);
\node[anchor=west,inner sep=0pt] at (14.85,17.5) {$v_{\xi^2}$};
\node[dot] at (g13) {};
\draw[black,line width=0.55pt] (g01) -- (7.45,16.68);
\draw[black,line width=0.55pt] (g01) -- (7.45,16.32);
\node[anchor=east,inner sep=0pt] at (7.15,16.5) {$u_{1}$};
\node[dot] at (g01) {};
\draw[black,line width=0.55pt] (g11) -- (14.55,16.68);
\draw[black,line width=0.55pt] (g11) -- (14.55,16.32);
\node[anchor=west,inner sep=0pt] at (14.85,16.5) {$v_{1}$};
\node[dot] at (g11) {};
\draw[black,line width=0.55pt] (g00) -- (7.45,15.68);
\draw[black,line width=0.55pt] (g00) -- (7.45,15.32);
\node[anchor=east,inner sep=0pt] at (7.15,15.5) {$u_{0}$};
\node[dot] at (g00) {};
\draw[black,line width=0.55pt] (g10) -- (14.55,15.68);
\draw[black,line width=0.55pt] (g10) -- (14.55,15.32);
\node[anchor=west,inner sep=0pt] at (14.85,15.5) {$v_{0}$};
\node[dot] at (g10) {};
\end{scope}
% Addition
\begin{scope}[shift={(8,15.4)}]
\foreach \i in {0,...,5} {
  \draw[line width=0.4pt] (1.5*\i,0) -- (1.5*\i,3.2);
  \draw[line width=0.4pt] (0,0.64*\i) -- (7.5,0.64*\i);
}
\node[inner sep=0pt] at (0.75,2.88) {$+$};
\node[inner sep=0pt] at (2.25,2.88) {$x+1$};
\node[inner sep=0pt] at (3.75,2.88) {$x$};
\node[inner sep=0pt] at (5.25,2.88) {$1$};
\node[inner sep=0pt] at (6.75,2.88) {$0$};
\node[inner sep=0pt] at (0.75,2.24) {$x+1$};
\node[inner sep=0pt] at (2.25,2.24) {$0$};
\node[inner sep=0pt] at (3.75,2.24) {$1$};
\node[inner sep=0pt] at (5.25,2.24) {$x$};
\node[inner sep=0pt] at (6.75,2.24) {$x+1$};
\node[inner sep=0pt] at (0.75,1.6) {$x$};
\node[inner sep=0pt] at (2.25,1.6) {$1$};
\node[inner sep=0pt] at (3.75,1.6) {$0$};
\node[inner sep=0pt] at (5.25,1.6) {$x+1$};
\node[inner sep=0pt] at (6.75,1.6) {$x$};
\node[inner sep=0pt] at (0.75,0.96) {$1$};
\node[inner sep=0pt] at (2.25,0.96) {$x$};
\node[inner sep=0pt] at (3.75,0.96) {$x+1$};
\node[inner sep=0pt] at (5.25,0.96) {$0$};
\node[inner sep=0pt] at (6.75,0.96) {$1$};
\node[inner sep=0pt] at (0.75,0.32) {$0$};
\node[inner sep=0pt] at (2.25,0.32) {$x+1$};
\node[inner sep=0pt] at (3.75,0.32) {$x$};
\node[inner sep=0pt] at (5.25,0.32) {$1$};
\node[inner sep=0pt] at (6.75,0.32) {$0$};
\end{scope}
\node at (12.5,14.5) {(b) Addition in $\mathbb{F}_4$};
% Multiplication
\begin{scope}[shift={(17.25,15.4)}]
\foreach \i in {0,...,5} {
  \draw[line width=0.4pt] (1.5*\i,0) -- (1.5*\i,3.2);
  \draw[line width=0.4pt] (0,0.64*\i) -- (7.5,0.64*\i);
}
\node[inner sep=0pt] at (0.75,2.88) {$\cdot$};
\node[inner sep=0pt] at (2.25,2.88) {$x+1$};
\node[inner sep=0pt] at (3.75,2.88) {$x$};
\node[inner sep=0pt] at (5.25,2.88) {$1$};
\node[inner sep=0pt] at (6.75,2.88) {$0$};
\node[inner sep=0pt] at (0.75,2.24) {$x+1$};
\node[inner sep=0pt] at (2.25,2.24) {$x$};
\node[inner sep=0pt] at (3.75,2.24) {$1$};
\node[inner sep=0pt] at (5.25,2.24) {$x+1$};
\node[inner sep=0pt] at (6.75,2.24) {$0$};
\node[inner sep=0pt] at (0.75,1.6) {$x$};
\node[inner sep=0pt] at (2.25,1.6) {$1$};
\node[inner sep=0pt] at (3.75,1.6) {$x+1$};
\node[inner sep=0pt] at (5.25,1.6) {$x$};
\node[inner sep=0pt] at (6.75,1.6) {$0$};
\node[inner sep=0pt] at (0.75,0.96) {$1$};
\node[inner sep=0pt] at (2.25,0.96) {$x+1$};
\node[inner sep=0pt] at (3.75,0.96) {$x$};
\node[inner sep=0pt] at (5.25,0.96) {$1$};
\node[inner sep=0pt] at (6.75,0.96) {$0$};
\node[inner sep=0pt] at (0.75,0.32) {$0$};
\node[inner sep=0pt] at (2.25,0.32) {$0$};
\node[inner sep=0pt] at (3.75,0.32) {$0$};
\node[inner sep=0pt] at (5.25,0.32) {$0$};
\node[inner sep=0pt] at (6.75,0.32) {$0$};
\end{scope}
\node at (21.5,14.5) {(c) Multiplication in $\mathbb{F}_4$};
\begin{scope}[shift={(0,-0)}]
\begin{scope}[yscale=0.85]
\node at (4,-1.2) {(d) The edges lifted from semiedges and}; \node at (4,-2) {~~~~~~from $M_1=\{u_{\xi}v_{\xi}, u_{\xi^2}v_{\xi^2}, u_{1}v_{1}, u_{0}v_{0}\}$};
\tikzset{dot/.style={circle,fill=black,draw=black,inner sep=0pt,minimum size=3pt}}
% In every fibre the order is xi, xi^2, 1, 0.
\coordinate (m1-h022) at (0,14.7);
\coordinate (m1-h023) at (0,13.9);
\coordinate (m1-h021) at (0,13.1);
\coordinate (m1-h020) at (0,12.3);
\coordinate (m1-h032) at (0,10.6);
\coordinate (m1-h033) at (0,9.8);
\coordinate (m1-h031) at (0,9);
\coordinate (m1-h030) at (0,8.2);
\coordinate (m1-h012) at (0,6.5);
\coordinate (m1-h013) at (0,5.7);
\coordinate (m1-h011) at (0,4.9);
\coordinate (m1-h010) at (0,4.1);
\coordinate (m1-h002) at (0,2.4);
\coordinate (m1-h003) at (0,1.6);
\coordinate (m1-h001) at (0,0.8);
\coordinate (m1-h000) at (0,0);
\coordinate (m1-h122) at (8,14.7);
\coordinate (m1-h123) at (8,13.9);
\coordinate (m1-h121) at (8,13.1);
\coordinate (m1-h120) at (8,12.3);
\coordinate (m1-h132) at (8,10.6);
\coordinate (m1-h133) at (8,9.8);
\coordinate (m1-h131) at (8,9);
\coordinate (m1-h130) at (8,8.2);
\coordinate (m1-h112) at (8,6.5);
\coordinate (m1-h113) at (8,5.7);
\coordinate (m1-h111) at (8,4.9);
\coordinate (m1-h110) at (8,4.1);
\coordinate (m1-h102) at (8,2.4);
\coordinate (m1-h103) at (8,1.6);
\coordinate (m1-h101) at (8,0.8);
\coordinate (m1-h100) at (8,0);
% The 16 cross edges of the lifted matching M_1.
\draw[draw=MMSblue,line width=0.65pt] (m1-h000) -- (m1-h100);
\draw[draw=MMSblue,line width=0.65pt] (m1-h001) -- (m1-h101);
\draw[draw=MMSblue,line width=0.65pt] (m1-h002) -- (m1-h102);
\draw[draw=MMSblue,line width=0.65pt] (m1-h003) -- (m1-h103);
\draw[draw=MMSblue,line width=0.65pt] (m1-h010) -- (m1-h111);
\draw[draw=MMSblue,line width=0.65pt] (m1-h011) -- (m1-h110);
\draw[draw=MMSblue,line width=0.65pt] (m1-h012) -- (m1-h113);
\draw[draw=MMSblue,line width=0.65pt] (m1-h013) -- (m1-h112);
\draw[draw=MMSblue,line width=0.65pt] (m1-h020) -- (m1-h123);
\draw[draw=MMSblue,line width=0.65pt] (m1-h021) -- (m1-h122);
\draw[draw=MMSblue,line width=0.65pt] (m1-h022) -- (m1-h121);
\draw[draw=MMSblue,line width=0.65pt] (m1-h023) -- (m1-h120);
\draw[draw=MMSblue,line width=0.65pt] (m1-h030) -- (m1-h132);
\draw[draw=MMSblue,line width=0.65pt] (m1-h031) -- (m1-h133);
\draw[draw=MMSblue,line width=0.65pt] (m1-h032) -- (m1-h130);
\draw[draw=MMSblue,line width=0.65pt] (m1-h033) -- (m1-h131);
% The 32 internal edges.
\draw[draw=black!45,line width=0.85pt] (m1-h022) .. controls (-0.17,14.46) and (-0.17,14.14) .. (m1-h023);
\draw[draw=black!45,line width=0.85pt] (m1-h021) .. controls (-0.17,12.86) and (-0.17,12.54) .. (m1-h020);
\draw[draw=black!45,line width=0.85pt] (m1-h022) .. controls (-0.7,14.332) and (-0.7,13.468) .. (m1-h021);
\draw[draw=black!45,line width=0.85pt,preaction={draw=white,line width=2.2pt}] (m1-h023) .. controls (-0.7,13.532) and (-0.7,12.668) .. (m1-h020);
\draw[draw=black!45,line width=0.85pt] (m1-h032) .. controls (-0.17,10.36) and (-0.17,10.04) .. (m1-h033);
\draw[draw=black!45,line width=0.85pt] (m1-h031) .. controls (-0.17,8.76) and (-0.17,8.44) .. (m1-h030);
\draw[draw=black!45,line width=0.85pt] (m1-h032) .. controls (-0.7,10.232) and (-0.7,9.368) .. (m1-h031);
\draw[draw=black!45,line width=0.85pt,preaction={draw=white,line width=2.2pt}] (m1-h033) .. controls (-0.7,9.432) and (-0.7,8.568) .. (m1-h030);
\draw[draw=black!45,line width=0.85pt] (m1-h012) .. controls (-0.17,6.26) and (-0.17,5.94) .. (m1-h013);
\draw[draw=black!45,line width=0.85pt] (m1-h011) .. controls (-0.17,4.66) and (-0.17,4.34) .. (m1-h010);
\draw[draw=black!45,line width=0.85pt] (m1-h012) .. controls (-0.7,6.132) and (-0.7,5.268) .. (m1-h011);
\draw[draw=black!45,line width=0.85pt,preaction={draw=white,line width=2.2pt}] (m1-h013) .. controls (-0.7,5.332) and (-0.7,4.468) .. (m1-h010);
\draw[draw=black!45,line width=0.85pt] (m1-h002) .. controls (-0.17,2.16) and (-0.17,1.84) .. (m1-h003);
\draw[draw=black!45,line width=0.85pt] (m1-h001) .. controls (-0.17,0.56) and (-0.17,0.24) .. (m1-h000);
\draw[draw=black!45,line width=0.85pt] (m1-h002) .. controls (-0.7,2.032) and (-0.7,1.168) .. (m1-h001);
\draw[draw=black!45,line width=0.85pt,preaction={draw=white,line width=2.2pt}] (m1-h003) .. controls (-0.7,1.232) and (-0.7,0.368) .. (m1-h000);
\draw[draw=black!45,line width=0.85pt] (m1-h122) .. controls (8.17,14.46) and (8.17,14.14) .. (m1-h123);
\draw[draw=black!45,line width=0.85pt] (m1-h123) .. controls (8.17,13.66) and (8.17,13.34) .. (m1-h121);
\draw[draw=black!45,line width=0.85pt] (m1-h121) .. controls (8.17,12.86) and (8.17,12.54) .. (m1-h120);
\draw[draw=black!45,line width=0.85pt] (m1-h122) .. controls (8.7,14.148) and (8.7,12.852) .. (m1-h120);
\draw[draw=black!45,line width=0.85pt] (m1-h132) .. controls (8.17,10.36) and (8.17,10.04) .. (m1-h133);
\draw[draw=black!45,line width=0.85pt] (m1-h133) .. controls (8.17,9.56) and (8.17,9.24) .. (m1-h131);
\draw[draw=black!45,line width=0.85pt] (m1-h131) .. controls (8.17,8.76) and (8.17,8.44) .. (m1-h130);
\draw[draw=black!45,line width=0.85pt] (m1-h132) .. controls (8.7,10.048) and (8.7,8.752) .. (m1-h130);
\draw[draw=black!45,line width=0.85pt] (m1-h112) .. controls (8.17,6.26) and (8.17,5.94) .. (m1-h113);
\draw[draw=black!45,line width=0.85pt] (m1-h113) .. controls (8.17,5.46) and (8.17,5.14) .. (m1-h111);
\draw[draw=black!45,line width=0.85pt] (m1-h111) .. controls (8.17,4.66) and (8.17,4.34) .. (m1-h110);
\draw[draw=black!45,line width=0.85pt] (m1-h112) .. controls (8.7,5.948) and (8.7,4.652) .. (m1-h110);
\draw[draw=black!45,line width=0.85pt] (m1-h102) .. controls (8.17,2.16) and (8.17,1.84) .. (m1-h103);
\draw[draw=black!45,line width=0.85pt] (m1-h103) .. controls (8.17,1.36) and (8.17,1.04) .. (m1-h101);
\draw[draw=black!45,line width=0.85pt] (m1-h101) .. controls (8.17,0.56) and (8.17,0.24) .. (m1-h100);
\draw[draw=black!45,line width=0.85pt] (m1-h102) .. controls (8.7,1.848) and (8.7,0.552) .. (m1-h100);
\node[dot] at (m1-h022) {};
\node[anchor=east,inner sep=0pt] at (-0.8,14.7) {$(u_{\xi},\xi)$};
\node[dot] at (m1-h023) {};
\node[anchor=east,inner sep=0pt] at (-0.8,13.9) {$(u_{\xi},\xi^2)$};
\node[dot] at (m1-h021) {};
\node[anchor=east,inner sep=0pt] at (-0.8,13.1) {$(u_{\xi},1)$};
\node[dot] at (m1-h020) {};
\node[anchor=east,inner sep=0pt] at (-0.8,12.3) {$(u_{\xi},0)$};
\node[dot] at (m1-h032) {};
\node[anchor=east,inner sep=0pt] at (-0.8,10.6) {$(u_{\xi^2},\xi)$};
\node[dot] at (m1-h033) {};
\node[anchor=east,inner sep=0pt] at (-0.8,9.8) {$(u_{\xi^2},\xi^2)$};
\node[dot] at (m1-h031) {};
\node[anchor=east,inner sep=0pt] at (-0.8,9) {$(u_{\xi^2},1)$};
\node[dot] at (m1-h030) {};
\node[anchor=east,inner sep=0pt] at (-0.8,8.2) {$(u_{\xi^2},0)$};
\node[dot] at (m1-h012) {};
\node[anchor=east,inner sep=0pt] at (-0.8,6.5) {$(u_{1},\xi)$};
\node[dot] at (m1-h013) {};
\node[anchor=east,inner sep=0pt] at (-0.8,5.7) {$(u_{1},\xi^2)$};
\node[dot] at (m1-h011) {};
\node[anchor=east,inner sep=0pt] at (-0.8,4.9) {$(u_{1},1)$};
\node[dot] at (m1-h010) {};
\node[anchor=east,inner sep=0pt] at (-0.8,4.1) {$(u_{1},0)$};
\node[dot] at (m1-h002) {};
\node[anchor=east,inner sep=0pt] at (-0.8,2.4) {$(u_{0},\xi)$};
\node[dot] at (m1-h003) {};
\node[anchor=east,inner sep=0pt] at (-0.8,1.6) {$(u_{0},\xi^2)$};
\node[dot] at (m1-h001) {};
\node[anchor=east,inner sep=0pt] at (-0.8,0.8) {$(u_{0},1)$};
\node[dot] at (m1-h000) {};
\node[anchor=east,inner sep=0pt] at (-0.8,0) {$(u_{0},0)$};
\node[dot] at (m1-h122) {};
\node[anchor=west,inner sep=0pt] at (8.8,14.7) {$(v_{\xi},\xi)$};
\node[dot] at (m1-h123) {};
\node[anchor=west,inner sep=0pt] at (8.8,13.9) {$(v_{\xi},\xi^2)$};
\node[dot] at (m1-h121) {};
\node[anchor=west,inner sep=0pt] at (8.8,13.1) {$(v_{\xi},1)$};
\node[dot] at (m1-h120) {};
\node[anchor=west,inner sep=0pt] at (8.8,12.3) {$(v_{\xi},0)$};
\node[dot] at (m1-h132) {};
\node[anchor=west,inner sep=0pt] at (8.8,10.6) {$(v_{\xi^2},\xi)$};
\node[dot] at (m1-h133) {};
\node[anchor=west,inner sep=0pt] at (8.8,9.8) {$(v_{\xi^2},\xi^2)$};
\node[dot] at (m1-h131) {};
\node[anchor=west,inner sep=0pt] at (8.8,9) {$(v_{\xi^2},1)$};
\node[dot] at (m1-h130) {};
\node[anchor=west,inner sep=0pt] at (8.8,8.2) {$(v_{\xi^2},0)$};
\node[dot] at (m1-h112) {};
\node[anchor=west,inner sep=0pt] at (8.8,6.5) {$(v_{1},\xi)$};
\node[dot] at (m1-h113) {};
\node[anchor=west,inner sep=0pt] at (8.8,5.7) {$(v_{1},\xi^2)$};
\node[dot] at (m1-h111) {};
\node[anchor=west,inner sep=0pt] at (8.8,4.9) {$(v_{1},1)$};
\node[dot] at (m1-h110) {};
\node[anchor=west,inner sep=0pt] at (8.8,4.1) {$(v_{1},0)$};
\node[dot] at (m1-h102) {};
\node[anchor=west,inner sep=0pt] at (8.8,2.4) {$(v_{0},\xi)$};
\node[dot] at (m1-h103) {};
\node[anchor=west,inner sep=0pt] at (8.8,1.6) {$(v_{0},\xi^2)$};
\node[dot] at (m1-h101) {};
\node[anchor=west,inner sep=0pt] at (8.8,0.8) {$(v_{0},1)$};
\node[dot] at (m1-h100) {};
\node[anchor=west,inner sep=0pt] at (8.8,0) {$(v_{0},0)$};
\end{scope}
\end{scope}
\begin{scope}[shift={(14,-0)}]
\begin{scope}[yscale=0.85]
\node at (4,-1.2) {(e) The edges lifted from semiedges and}; \node at (4,-2) {~~~~~~from $M_2=\{u_{\xi}v_{\xi^2}, u_{\xi^2}v_{1}, u_{1}v_{0}, u_{0}v_{\xi}\}$};
\tikzset{dot/.style={circle,fill=black,draw=black,inner sep=0pt,minimum size=3pt}}
% In every fibre the order is xi, xi^2, 1, 0.
\coordinate (m2-h022) at (0,14.7);
\coordinate (m2-h023) at (0,13.9);
\coordinate (m2-h021) at (0,13.1);
\coordinate (m2-h020) at (0,12.3);
\coordinate (m2-h032) at (0,10.6);
\coordinate (m2-h033) at (0,9.8);
\coordinate (m2-h031) at (0,9);
\coordinate (m2-h030) at (0,8.2);
\coordinate (m2-h012) at (0,6.5);
\coordinate (m2-h013) at (0,5.7);
\coordinate (m2-h011) at (0,4.9);
\coordinate (m2-h010) at (0,4.1);
\coordinate (m2-h002) at (0,2.4);
\coordinate (m2-h003) at (0,1.6);
\coordinate (m2-h001) at (0,0.8);
\coordinate (m2-h000) at (0,0);
\coordinate (m2-h122) at (8,14.7);
\coordinate (m2-h123) at (8,13.9);
\coordinate (m2-h121) at (8,13.1);
\coordinate (m2-h120) at (8,12.3);
\coordinate (m2-h132) at (8,10.6);
\coordinate (m2-h133) at (8,9.8);
\coordinate (m2-h131) at (8,9);
\coordinate (m2-h130) at (8,8.2);
\coordinate (m2-h112) at (8,6.5);
\coordinate (m2-h113) at (8,5.7);
\coordinate (m2-h111) at (8,4.9);
\coordinate (m2-h110) at (8,4.1);
\coordinate (m2-h102) at (8,2.4);
\coordinate (m2-h103) at (8,1.6);
\coordinate (m2-h101) at (8,0.8);
\coordinate (m2-h100) at (8,0);
% The 16 cross edges of the lifted matching M_2.
\draw[draw=MMSorange,line width=0.65pt] (m2-h000) -- (m2-h120);
\draw[draw=MMSorange,line width=0.65pt] (m2-h001) -- (m2-h121);
\draw[draw=MMSorange,line width=0.65pt] (m2-h002) -- (m2-h122);
\draw[draw=MMSorange,line width=0.65pt] (m2-h003) -- (m2-h123);
\draw[draw=MMSorange,line width=0.65pt] (m2-h010) -- (m2-h100);
\draw[draw=MMSorange,line width=0.65pt] (m2-h011) -- (m2-h101);
\draw[draw=MMSorange,line width=0.65pt] (m2-h012) -- (m2-h102);
\draw[draw=MMSorange,line width=0.65pt] (m2-h013) -- (m2-h103);
\draw[draw=MMSorange,line width=0.65pt] (m2-h020) -- (m2-h131);
\draw[draw=MMSorange,line width=0.65pt] (m2-h021) -- (m2-h130);
\draw[draw=MMSorange,line width=0.65pt] (m2-h022) -- (m2-h133);
\draw[draw=MMSorange,line width=0.65pt] (m2-h023) -- (m2-h132);
\draw[draw=MMSorange,line width=0.65pt] (m2-h030) -- (m2-h113);
\draw[draw=MMSorange,line width=0.65pt] (m2-h031) -- (m2-h112);
\draw[draw=MMSorange,line width=0.65pt] (m2-h032) -- (m2-h111);
\draw[draw=MMSorange,line width=0.65pt] (m2-h033) -- (m2-h110);
% The 32 internal edges.
\draw[draw=black!45,line width=0.85pt] (m2-h022) .. controls (-0.17,14.46) and (-0.17,14.14) .. (m2-h023);
\draw[draw=black!45,line width=0.85pt] (m2-h021) .. controls (-0.17,12.86) and (-0.17,12.54) .. (m2-h020);
\draw[draw=black!45,line width=0.85pt] (m2-h022) .. controls (-0.7,14.332) and (-0.7,13.468) .. (m2-h021);
\draw[draw=black!45,line width=0.85pt,preaction={draw=white,line width=2.2pt}] (m2-h023) .. controls (-0.7,13.532) and (-0.7,12.668) .. (m2-h020);
\draw[draw=black!45,line width=0.85pt] (m2-h032) .. controls (-0.17,10.36) and (-0.17,10.04) .. (m2-h033);
\draw[draw=black!45,line width=0.85pt] (m2-h031) .. controls (-0.17,8.76) and (-0.17,8.44) .. (m2-h030);
\draw[draw=black!45,line width=0.85pt] (m2-h032) .. controls (-0.7,10.232) and (-0.7,9.368) .. (m2-h031);
\draw[draw=black!45,line width=0.85pt,preaction={draw=white,line width=2.2pt}] (m2-h033) .. controls (-0.7,9.432) and (-0.7,8.568) .. (m2-h030);
\draw[draw=black!45,line width=0.85pt] (m2-h012) .. controls (-0.17,6.26) and (-0.17,5.94) .. (m2-h013);
\draw[draw=black!45,line width=0.85pt] (m2-h011) .. controls (-0.17,4.66) and (-0.17,4.34) .. (m2-h010);
\draw[draw=black!45,line width=0.85pt] (m2-h012) .. controls (-0.7,6.132) and (-0.7,5.268) .. (m2-h011);
\draw[draw=black!45,line width=0.85pt,preaction={draw=white,line width=2.2pt}] (m2-h013) .. controls (-0.7,5.332) and (-0.7,4.468) .. (m2-h010);
\draw[draw=black!45,line width=0.85pt] (m2-h002) .. controls (-0.17,2.16) and (-0.17,1.84) .. (m2-h003);
\draw[draw=black!45,line width=0.85pt] (m2-h001) .. controls (-0.17,0.56) and (-0.17,0.24) .. (m2-h000);
\draw[draw=black!45,line width=0.85pt] (m2-h002) .. controls (-0.7,2.032) and (-0.7,1.168) .. (m2-h001);
\draw[draw=black!45,line width=0.85pt,preaction={draw=white,line width=2.2pt}] (m2-h003) .. controls (-0.7,1.232) and (-0.7,0.368) .. (m2-h000);
\draw[draw=black!45,line width=0.85pt] (m2-h122) .. controls (8.17,14.46) and (8.17,14.14) .. (m2-h123);
\draw[draw=black!45,line width=0.85pt] (m2-h123) .. controls (8.17,13.66) and (8.17,13.34) .. (m2-h121);
\draw[draw=black!45,line width=0.85pt] (m2-h121) .. controls (8.17,12.86) and (8.17,12.54) .. (m2-h120);
\draw[draw=black!45,line width=0.85pt] (m2-h122) .. controls (8.7,14.148) and (8.7,12.852) .. (m2-h120);
\draw[draw=black!45,line width=0.85pt] (m2-h132) .. controls (8.17,10.36) and (8.17,10.04) .. (m2-h133);
\draw[draw=black!45,line width=0.85pt] (m2-h133) .. controls (8.17,9.56) and (8.17,9.24) .. (m2-h131);
\draw[draw=black!45,line width=0.85pt] (m2-h131) .. controls (8.17,8.76) and (8.17,8.44) .. (m2-h130);
\draw[draw=black!45,line width=0.85pt] (m2-h132) .. controls (8.7,10.048) and (8.7,8.752) .. (m2-h130);
\draw[draw=black!45,line width=0.85pt] (m2-h112) .. controls (8.17,6.26) and (8.17,5.94) .. (m2-h113);
\draw[draw=black!45,line width=0.85pt] (m2-h113) .. controls (8.17,5.46) and (8.17,5.14) .. (m2-h111);
\draw[draw=black!45,line width=0.85pt] (m2-h111) .. controls (8.17,4.66) and (8.17,4.34) .. (m2-h110);
\draw[draw=black!45,line width=0.85pt] (m2-h112) .. controls (8.7,5.948) and (8.7,4.652) .. (m2-h110);
\draw[draw=black!45,line width=0.85pt] (m2-h102) .. controls (8.17,2.16) and (8.17,1.84) .. (m2-h103);
\draw[draw=black!45,line width=0.85pt] (m2-h103) .. controls (8.17,1.36) and (8.17,1.04) .. (m2-h101);
\draw[draw=black!45,line width=0.85pt] (m2-h101) .. controls (8.17,0.56) and (8.17,0.24) .. (m2-h100);
\draw[draw=black!45,line width=0.85pt] (m2-h102) .. controls (8.7,1.848) and (8.7,0.552) .. (m2-h100);
\node[dot] at (m2-h022) {};
\node[anchor=east,inner sep=0pt] at (-0.8,14.7) {$(u_{\xi},\xi)$};
\node[dot] at (m2-h023) {};
\node[anchor=east,inner sep=0pt] at (-0.8,13.9) {$(u_{\xi},\xi^2)$};
\node[dot] at (m2-h021) {};
\node[anchor=east,inner sep=0pt] at (-0.8,13.1) {$(u_{\xi},1)$};
\node[dot] at (m2-h020) {};
\node[anchor=east,inner sep=0pt] at (-0.8,12.3) {$(u_{\xi},0)$};
\node[dot] at (m2-h032) {};
\node[anchor=east,inner sep=0pt] at (-0.8,10.6) {$(u_{\xi^2},\xi)$};
\node[dot] at (m2-h033) {};
\node[anchor=east,inner sep=0pt] at (-0.8,9.8) {$(u_{\xi^2},\xi^2)$};
\node[dot] at (m2-h031) {};
\node[anchor=east,inner sep=0pt] at (-0.8,9) {$(u_{\xi^2},1)$};
\node[dot] at (m2-h030) {};
\node[anchor=east,inner sep=0pt] at (-0.8,8.2) {$(u_{\xi^2},0)$};
\node[dot] at (m2-h012) {};
\node[anchor=east,inner sep=0pt] at (-0.8,6.5) {$(u_{1},\xi)$};
\node[dot] at (m2-h013) {};
\node[anchor=east,inner sep=0pt] at (-0.8,5.7) {$(u_{1},\xi^2)$};
\node[dot] at (m2-h011) {};
\node[anchor=east,inner sep=0pt] at (-0.8,4.9) {$(u_{1},1)$};
\node[dot] at (m2-h010) {};
\node[anchor=east,inner sep=0pt] at (-0.8,4.1) {$(u_{1},0)$};
\node[dot] at (m2-h002) {};
\node[anchor=east,inner sep=0pt] at (-0.8,2.4) {$(u_{0},\xi)$};
\node[dot] at (m2-h003) {};
\node[anchor=east,inner sep=0pt] at (-0.8,1.6) {$(u_{0},\xi^2)$};
\node[dot] at (m2-h001) {};
\node[anchor=east,inner sep=0pt] at (-0.8,0.8) {$(u_{0},1)$};
\node[dot] at (m2-h000) {};
\node[anchor=east,inner sep=0pt] at (-0.8,0) {$(u_{0},0)$};
\node[dot] at (m2-h122) {};
\node[anchor=west,inner sep=0pt] at (8.8,14.7) {$(v_{\xi},\xi)$};
\node[dot] at (m2-h123) {};
\node[anchor=west,inner sep=0pt] at (8.8,13.9) {$(v_{\xi},\xi^2)$};
\node[dot] at (m2-h121) {};
\node[anchor=west,inner sep=0pt] at (8.8,13.1) {$(v_{\xi},1)$};
\node[dot] at (m2-h120) {};
\node[anchor=west,inner sep=0pt] at (8.8,12.3) {$(v_{\xi},0)$};
\node[dot] at (m2-h132) {};
\node[anchor=west,inner sep=0pt] at (8.8,10.6) {$(v_{\xi^2},\xi)$};
\node[dot] at (m2-h133) {};
\node[anchor=west,inner sep=0pt] at (8.8,9.8) {$(v_{\xi^2},\xi^2)$};
\node[dot] at (m2-h131) {};
\node[anchor=west,inner sep=0pt] at (8.8,9) {$(v_{\xi^2},1)$};
\node[dot] at (m2-h130) {};
\node[anchor=west,inner sep=0pt] at (8.8,8.2) {$(v_{\xi^2},0)$};
\node[dot] at (m2-h112) {};
\node[anchor=west,inner sep=0pt] at (8.8,6.5) {$(v_{1},\xi)$};
\node[dot] at (m2-h113) {};
\node[anchor=west,inner sep=0pt] at (8.8,5.7) {$(v_{1},\xi^2)$};
\node[dot] at (m2-h111) {};
\node[anchor=west,inner sep=0pt] at (8.8,4.9) {$(v_{1},1)$};
\node[dot] at (m2-h110) {};
\node[anchor=west,inner sep=0pt] at (8.8,4.1) {$(v_{1},0)$};
\node[dot] at (m2-h102) {};
\node[anchor=west,inner sep=0pt] at (8.8,2.4) {$(v_{0},\xi)$};
\node[dot] at (m2-h103) {};
\node[anchor=west,inner sep=0pt] at (8.8,1.6) {$(v_{0},\xi^2)$};
\node[dot] at (m2-h101) {};
\node[anchor=west,inner sep=0pt] at (8.8,0.8) {$(v_{0},1)$};
\node[dot] at (m2-h100) {};
\node[anchor=west,inner sep=0pt] at (8.8,0) {$(v_{0},0)$};
\end{scope}
\end{scope}
\begin{scope}[shift={(0,-15.5)}]
\begin{scope}[yscale=0.85]
\node at (4,-1.2) {(f) The edges lifted from semiedges and}; \node at (4,-2) {~~~~~~from $M_3=\{u_{\xi}v_{1}, u_{\xi^2}v_{0}, u_{1}v_{\xi}, u_{0}v_{\xi^2}\}$};
\tikzset{dot/.style={circle,fill=black,draw=black,inner sep=0pt,minimum size=3pt}}
% In every fibre the order is xi, xi^2, 1, 0.
\coordinate (m3-h022) at (0,14.7);
\coordinate (m3-h023) at (0,13.9);
\coordinate (m3-h021) at (0,13.1);
\coordinate (m3-h020) at (0,12.3);
\coordinate (m3-h032) at (0,10.6);
\coordinate (m3-h033) at (0,9.8);
\coordinate (m3-h031) at (0,9);
\coordinate (m3-h030) at (0,8.2);
\coordinate (m3-h012) at (0,6.5);
\coordinate (m3-h013) at (0,5.7);
\coordinate (m3-h011) at (0,4.9);
\coordinate (m3-h010) at (0,4.1);
\coordinate (m3-h002) at (0,2.4);
\coordinate (m3-h003) at (0,1.6);
\coordinate (m3-h001) at (0,0.8);
\coordinate (m3-h000) at (0,0);
\coordinate (m3-h122) at (8,14.7);
\coordinate (m3-h123) at (8,13.9);
\coordinate (m3-h121) at (8,13.1);
\coordinate (m3-h120) at (8,12.3);
\coordinate (m3-h132) at (8,10.6);
\coordinate (m3-h133) at (8,9.8);
\coordinate (m3-h131) at (8,9);
\coordinate (m3-h130) at (8,8.2);
\coordinate (m3-h112) at (8,6.5);
\coordinate (m3-h113) at (8,5.7);
\coordinate (m3-h111) at (8,4.9);
\coordinate (m3-h110) at (8,4.1);
\coordinate (m3-h102) at (8,2.4);
\coordinate (m3-h103) at (8,1.6);
\coordinate (m3-h101) at (8,0.8);
\coordinate (m3-h100) at (8,0);
% The 16 cross edges of the lifted matching M_3.
\draw[draw=MMSgreen,line width=0.65pt] (m3-h000) -- (m3-h130);
\draw[draw=MMSgreen,line width=0.65pt] (m3-h001) -- (m3-h131);
\draw[draw=MMSgreen,line width=0.65pt] (m3-h002) -- (m3-h132);
\draw[draw=MMSgreen,line width=0.65pt] (m3-h003) -- (m3-h133);
\draw[draw=MMSgreen,line width=0.65pt] (m3-h010) -- (m3-h122);
\draw[draw=MMSgreen,line width=0.65pt] (m3-h011) -- (m3-h123);
\draw[draw=MMSgreen,line width=0.65pt] (m3-h012) -- (m3-h120);
\draw[draw=MMSgreen,line width=0.65pt] (m3-h013) -- (m3-h121);
\draw[draw=MMSgreen,line width=0.65pt] (m3-h020) -- (m3-h112);
\draw[draw=MMSgreen,line width=0.65pt] (m3-h021) -- (m3-h113);
\draw[draw=MMSgreen,line width=0.65pt] (m3-h022) -- (m3-h110);
\draw[draw=MMSgreen,line width=0.65pt] (m3-h023) -- (m3-h111);
\draw[draw=MMSgreen,line width=0.65pt] (m3-h030) -- (m3-h100);
\draw[draw=MMSgreen,line width=0.65pt] (m3-h031) -- (m3-h101);
\draw[draw=MMSgreen,line width=0.65pt] (m3-h032) -- (m3-h102);
\draw[draw=MMSgreen,line width=0.65pt] (m3-h033) -- (m3-h103);
% The 32 internal edges.
\draw[draw=black!45,line width=0.85pt] (m3-h022) .. controls (-0.17,14.46) and (-0.17,14.14) .. (m3-h023);
\draw[draw=black!45,line width=0.85pt] (m3-h021) .. controls (-0.17,12.86) and (-0.17,12.54) .. (m3-h020);
\draw[draw=black!45,line width=0.85pt] (m3-h022) .. controls (-0.7,14.332) and (-0.7,13.468) .. (m3-h021);
\draw[draw=black!45,line width=0.85pt,preaction={draw=white,line width=2.2pt}] (m3-h023) .. controls (-0.7,13.532) and (-0.7,12.668) .. (m3-h020);
\draw[draw=black!45,line width=0.85pt] (m3-h032) .. controls (-0.17,10.36) and (-0.17,10.04) .. (m3-h033);
\draw[draw=black!45,line width=0.85pt] (m3-h031) .. controls (-0.17,8.76) and (-0.17,8.44) .. (m3-h030);
\draw[draw=black!45,line width=0.85pt] (m3-h032) .. controls (-0.7,10.232) and (-0.7,9.368) .. (m3-h031);
\draw[draw=black!45,line width=0.85pt,preaction={draw=white,line width=2.2pt}] (m3-h033) .. controls (-0.7,9.432) and (-0.7,8.568) .. (m3-h030);
\draw[draw=black!45,line width=0.85pt] (m3-h012) .. controls (-0.17,6.26) and (-0.17,5.94) .. (m3-h013);
\draw[draw=black!45,line width=0.85pt] (m3-h011) .. controls (-0.17,4.66) and (-0.17,4.34) .. (m3-h010);
\draw[draw=black!45,line width=0.85pt] (m3-h012) .. controls (-0.7,6.132) and (-0.7,5.268) .. (m3-h011);
\draw[draw=black!45,line width=0.85pt,preaction={draw=white,line width=2.2pt}] (m3-h013) .. controls (-0.7,5.332) and (-0.7,4.468) .. (m3-h010);
\draw[draw=black!45,line width=0.85pt] (m3-h002) .. controls (-0.17,2.16) and (-0.17,1.84) .. (m3-h003);
\draw[draw=black!45,line width=0.85pt] (m3-h001) .. controls (-0.17,0.56) and (-0.17,0.24) .. (m3-h000);
\draw[draw=black!45,line width=0.85pt] (m3-h002) .. controls (-0.7,2.032) and (-0.7,1.168) .. (m3-h001);
\draw[draw=black!45,line width=0.85pt,preaction={draw=white,line width=2.2pt}] (m3-h003) .. controls (-0.7,1.232) and (-0.7,0.368) .. (m3-h000);
\draw[draw=black!45,line width=0.85pt] (m3-h122) .. controls (8.17,14.46) and (8.17,14.14) .. (m3-h123);
\draw[draw=black!45,line width=0.85pt] (m3-h123) .. controls (8.17,13.66) and (8.17,13.34) .. (m3-h121);
\draw[draw=black!45,line width=0.85pt] (m3-h121) .. controls (8.17,12.86) and (8.17,12.54) .. (m3-h120);
\draw[draw=black!45,line width=0.85pt] (m3-h122) .. controls (8.7,14.148) and (8.7,12.852) .. (m3-h120);
\draw[draw=black!45,line width=0.85pt] (m3-h132) .. controls (8.17,10.36) and (8.17,10.04) .. (m3-h133);
\draw[draw=black!45,line width=0.85pt] (m3-h133) .. controls (8.17,9.56) and (8.17,9.24) .. (m3-h131);
\draw[draw=black!45,line width=0.85pt] (m3-h131) .. controls (8.17,8.76) and (8.17,8.44) .. (m3-h130);
\draw[draw=black!45,line width=0.85pt] (m3-h132) .. controls (8.7,10.048) and (8.7,8.752) .. (m3-h130);
\draw[draw=black!45,line width=0.85pt] (m3-h112) .. controls (8.17,6.26) and (8.17,5.94) .. (m3-h113);
\draw[draw=black!45,line width=0.85pt] (m3-h113) .. controls (8.17,5.46) and (8.17,5.14) .. (m3-h111);
\draw[draw=black!45,line width=0.85pt] (m3-h111) .. controls (8.17,4.66) and (8.17,4.34) .. (m3-h110);
\draw[draw=black!45,line width=0.85pt] (m3-h112) .. controls (8.7,5.948) and (8.7,4.652) .. (m3-h110);
\draw[draw=black!45,line width=0.85pt] (m3-h102) .. controls (8.17,2.16) and (8.17,1.84) .. (m3-h103);
\draw[draw=black!45,line width=0.85pt] (m3-h103) .. controls (8.17,1.36) and (8.17,1.04) .. (m3-h101);
\draw[draw=black!45,line width=0.85pt] (m3-h101) .. controls (8.17,0.56) and (8.17,0.24) .. (m3-h100);
\draw[draw=black!45,line width=0.85pt] (m3-h102) .. controls (8.7,1.848) and (8.7,0.552) .. (m3-h100);
\node[dot] at (m3-h022) {};
\node[anchor=east,inner sep=0pt] at (-0.8,14.7) {$(u_{\xi},\xi)$};
\node[dot] at (m3-h023) {};
\node[anchor=east,inner sep=0pt] at (-0.8,13.9) {$(u_{\xi},\xi^2)$};
\node[dot] at (m3-h021) {};
\node[anchor=east,inner sep=0pt] at (-0.8,13.1) {$(u_{\xi},1)$};
\node[dot] at (m3-h020) {};
\node[anchor=east,inner sep=0pt] at (-0.8,12.3) {$(u_{\xi},0)$};
\node[dot] at (m3-h032) {};
\node[anchor=east,inner sep=0pt] at (-0.8,10.6) {$(u_{\xi^2},\xi)$};
\node[dot] at (m3-h033) {};
\node[anchor=east,inner sep=0pt] at (-0.8,9.8) {$(u_{\xi^2},\xi^2)$};
\node[dot] at (m3-h031) {};
\node[anchor=east,inner sep=0pt] at (-0.8,9) {$(u_{\xi^2},1)$};
\node[dot] at (m3-h030) {};
\node[anchor=east,inner sep=0pt] at (-0.8,8.2) {$(u_{\xi^2},0)$};
\node[dot] at (m3-h012) {};
\node[anchor=east,inner sep=0pt] at (-0.8,6.5) {$(u_{1},\xi)$};
\node[dot] at (m3-h013) {};
\node[anchor=east,inner sep=0pt] at (-0.8,5.7) {$(u_{1},\xi^2)$};
\node[dot] at (m3-h011) {};
\node[anchor=east,inner sep=0pt] at (-0.8,4.9) {$(u_{1},1)$};
\node[dot] at (m3-h010) {};
\node[anchor=east,inner sep=0pt] at (-0.8,4.1) {$(u_{1},0)$};
\node[dot] at (m3-h002) {};
\node[anchor=east,inner sep=0pt] at (-0.8,2.4) {$(u_{0},\xi)$};
\node[dot] at (m3-h003) {};
\node[anchor=east,inner sep=0pt] at (-0.8,1.6) {$(u_{0},\xi^2)$};
\node[dot] at (m3-h001) {};
\node[anchor=east,inner sep=0pt] at (-0.8,0.8) {$(u_{0},1)$};
\node[dot] at (m3-h000) {};
\node[anchor=east,inner sep=0pt] at (-0.8,0) {$(u_{0},0)$};
\node[dot] at (m3-h122) {};
\node[anchor=west,inner sep=0pt] at (8.8,14.7) {$(v_{\xi},\xi)$};
\node[dot] at (m3-h123) {};
\node[anchor=west,inner sep=0pt] at (8.8,13.9) {$(v_{\xi},\xi^2)$};
\node[dot] at (m3-h121) {};
\node[anchor=west,inner sep=0pt] at (8.8,13.1) {$(v_{\xi},1)$};
\node[dot] at (m3-h120) {};
\node[anchor=west,inner sep=0pt] at (8.8,12.3) {$(v_{\xi},0)$};
\node[dot] at (m3-h132) {};
\node[anchor=west,inner sep=0pt] at (8.8,10.6) {$(v_{\xi^2},\xi)$};
\node[dot] at (m3-h133) {};
\node[anchor=west,inner sep=0pt] at (8.8,9.8) {$(v_{\xi^2},\xi^2)$};
\node[dot] at (m3-h131) {};
\node[anchor=west,inner sep=0pt] at (8.8,9) {$(v_{\xi^2},1)$};
\node[dot] at (m3-h130) {};
\node[anchor=west,inner sep=0pt] at (8.8,8.2) {$(v_{\xi^2},0)$};
\node[dot] at (m3-h112) {};
\node[anchor=west,inner sep=0pt] at (8.8,6.5) {$(v_{1},\xi)$};
\node[dot] at (m3-h113) {};
\node[anchor=west,inner sep=0pt] at (8.8,5.7) {$(v_{1},\xi^2)$};
\node[dot] at (m3-h111) {};
\node[anchor=west,inner sep=0pt] at (8.8,4.9) {$(v_{1},1)$};
\node[dot] at (m3-h110) {};
\node[anchor=west,inner sep=0pt] at (8.8,4.1) {$(v_{1},0)$};
\node[dot] at (m3-h102) {};
\node[anchor=west,inner sep=0pt] at (8.8,2.4) {$(v_{0},\xi)$};
\node[dot] at (m3-h103) {};
\node[anchor=west,inner sep=0pt] at (8.8,1.6) {$(v_{0},\xi^2)$};
\node[dot] at (m3-h101) {};
\node[anchor=west,inner sep=0pt] at (8.8,0.8) {$(v_{0},1)$};
\node[dot] at (m3-h100) {};
\node[anchor=west,inner sep=0pt] at (8.8,0) {$(v_{0},0)$};
\end{scope}
\end{scope}
\begin{scope}[shift={(14,-15.5)}]
\begin{scope}[yscale=0.85]
\node at (4,-1.2) {(g) The edges lifted from semiedges and}; \node at (4,-2) {~~~~~~from $M_4=\{u_{\xi}v_{0}, u_{\xi^2}v_{\xi}, u_{1}v_{\xi^2}, u_{0}v_{1}\}$};
\tikzset{dot/.style={circle,fill=black,draw=black,inner sep=0pt,minimum size=3pt}}
% In every fibre the order is xi, xi^2, 1, 0.
\coordinate (m4-h022) at (0,14.7);
\coordinate (m4-h023) at (0,13.9);
\coordinate (m4-h021) at (0,13.1);
\coordinate (m4-h020) at (0,12.3);
\coordinate (m4-h032) at (0,10.6);
\coordinate (m4-h033) at (0,9.8);
\coordinate (m4-h031) at (0,9);
\coordinate (m4-h030) at (0,8.2);
\coordinate (m4-h012) at (0,6.5);
\coordinate (m4-h013) at (0,5.7);
\coordinate (m4-h011) at (0,4.9);
\coordinate (m4-h010) at (0,4.1);
\coordinate (m4-h002) at (0,2.4);
\coordinate (m4-h003) at (0,1.6);
\coordinate (m4-h001) at (0,0.8);
\coordinate (m4-h000) at (0,0);
\coordinate (m4-h122) at (8,14.7);
\coordinate (m4-h123) at (8,13.9);
\coordinate (m4-h121) at (8,13.1);
\coordinate (m4-h120) at (8,12.3);
\coordinate (m4-h132) at (8,10.6);
\coordinate (m4-h133) at (8,9.8);
\coordinate (m4-h131) at (8,9);
\coordinate (m4-h130) at (8,8.2);
\coordinate (m4-h112) at (8,6.5);
\coordinate (m4-h113) at (8,5.7);
\coordinate (m4-h111) at (8,4.9);
\coordinate (m4-h110) at (8,4.1);
\coordinate (m4-h102) at (8,2.4);
\coordinate (m4-h103) at (8,1.6);
\coordinate (m4-h101) at (8,0.8);
\coordinate (m4-h100) at (8,0);
% The 16 cross edges of the lifted matching M_4.
\draw[draw=MMSpurple,line width=0.65pt] (m4-h000) -- (m4-h110);
\draw[draw=MMSpurple,line width=0.65pt] (m4-h001) -- (m4-h111);
\draw[draw=MMSpurple,line width=0.65pt] (m4-h002) -- (m4-h112);
\draw[draw=MMSpurple,line width=0.65pt] (m4-h003) -- (m4-h113);
\draw[draw=MMSpurple,line width=0.65pt] (m4-h010) -- (m4-h133);
\draw[draw=MMSpurple,line width=0.65pt] (m4-h011) -- (m4-h132);
\draw[draw=MMSpurple,line width=0.65pt] (m4-h012) -- (m4-h131);
\draw[draw=MMSpurple,line width=0.65pt] (m4-h013) -- (m4-h130);
\draw[draw=MMSpurple,line width=0.65pt] (m4-h020) -- (m4-h100);
\draw[draw=MMSpurple,line width=0.65pt] (m4-h021) -- (m4-h101);
\draw[draw=MMSpurple,line width=0.65pt] (m4-h022) -- (m4-h102);
\draw[draw=MMSpurple,line width=0.65pt] (m4-h023) -- (m4-h103);
\draw[draw=MMSpurple,line width=0.65pt] (m4-h030) -- (m4-h121);
\draw[draw=MMSpurple,line width=0.65pt] (m4-h031) -- (m4-h120);
\draw[draw=MMSpurple,line width=0.65pt] (m4-h032) -- (m4-h123);
\draw[draw=MMSpurple,line width=0.65pt] (m4-h033) -- (m4-h122);
% The 32 internal edges.
\draw[draw=black!45,line width=0.85pt] (m4-h022) .. controls (-0.17,14.46) and (-0.17,14.14) .. (m4-h023);
\draw[draw=black!45,line width=0.85pt] (m4-h021) .. controls (-0.17,12.86) and (-0.17,12.54) .. (m4-h020);
\draw[draw=black!45,line width=0.85pt] (m4-h022) .. controls (-0.7,14.332) and (-0.7,13.468) .. (m4-h021);
\draw[draw=black!45,line width=0.85pt,preaction={draw=white,line width=2.2pt}] (m4-h023) .. controls (-0.7,13.532) and (-0.7,12.668) .. (m4-h020);
\draw[draw=black!45,line width=0.85pt] (m4-h032) .. controls (-0.17,10.36) and (-0.17,10.04) .. (m4-h033);
\draw[draw=black!45,line width=0.85pt] (m4-h031) .. controls (-0.17,8.76) and (-0.17,8.44) .. (m4-h030);
\draw[draw=black!45,line width=0.85pt] (m4-h032) .. controls (-0.7,10.232) and (-0.7,9.368) .. (m4-h031);
\draw[draw=black!45,line width=0.85pt,preaction={draw=white,line width=2.2pt}] (m4-h033) .. controls (-0.7,9.432) and (-0.7,8.568) .. (m4-h030);
\draw[draw=black!45,line width=0.85pt] (m4-h012) .. controls (-0.17,6.26) and (-0.17,5.94) .. (m4-h013);
\draw[draw=black!45,line width=0.85pt] (m4-h011) .. controls (-0.17,4.66) and (-0.17,4.34) .. (m4-h010);
\draw[draw=black!45,line width=0.85pt] (m4-h012) .. controls (-0.7,6.132) and (-0.7,5.268) .. (m4-h011);
\draw[draw=black!45,line width=0.85pt,preaction={draw=white,line width=2.2pt}] (m4-h013) .. controls (-0.7,5.332) and (-0.7,4.468) .. (m4-h010);
\draw[draw=black!45,line width=0.85pt] (m4-h002) .. controls (-0.17,2.16) and (-0.17,1.84) .. (m4-h003);
\draw[draw=black!45,line width=0.85pt] (m4-h001) .. controls (-0.17,0.56) and (-0.17,0.24) .. (m4-h000);
\draw[draw=black!45,line width=0.85pt] (m4-h002) .. controls (-0.7,2.032) and (-0.7,1.168) .. (m4-h001);
\draw[draw=black!45,line width=0.85pt,preaction={draw=white,line width=2.2pt}] (m4-h003) .. controls (-0.7,1.232) and (-0.7,0.368) .. (m4-h000);
\draw[draw=black!45,line width=0.85pt] (m4-h122) .. controls (8.17,14.46) and (8.17,14.14) .. (m4-h123);
\draw[draw=black!45,line width=0.85pt] (m4-h123) .. controls (8.17,13.66) and (8.17,13.34) .. (m4-h121);
\draw[draw=black!45,line width=0.85pt] (m4-h121) .. controls (8.17,12.86) and (8.17,12.54) .. (m4-h120);
\draw[draw=black!45,line width=0.85pt] (m4-h122) .. controls (8.7,14.148) and (8.7,12.852) .. (m4-h120);
\draw[draw=black!45,line width=0.85pt] (m4-h132) .. controls (8.17,10.36) and (8.17,10.04) .. (m4-h133);
\draw[draw=black!45,line width=0.85pt] (m4-h133) .. controls (8.17,9.56) and (8.17,9.24) .. (m4-h131);
\draw[draw=black!45,line width=0.85pt] (m4-h131) .. controls (8.17,8.76) and (8.17,8.44) .. (m4-h130);
\draw[draw=black!45,line width=0.85pt] (m4-h132) .. controls (8.7,10.048) and (8.7,8.752) .. (m4-h130);
\draw[draw=black!45,line width=0.85pt] (m4-h112) .. controls (8.17,6.26) and (8.17,5.94) .. (m4-h113);
\draw[draw=black!45,line width=0.85pt] (m4-h113) .. controls (8.17,5.46) and (8.17,5.14) .. (m4-h111);
\draw[draw=black!45,line width=0.85pt] (m4-h111) .. controls (8.17,4.66) and (8.17,4.34) .. (m4-h110);
\draw[draw=black!45,line width=0.85pt] (m4-h112) .. controls (8.7,5.948) and (8.7,4.652) .. (m4-h110);
\draw[draw=black!45,line width=0.85pt] (m4-h102) .. controls (8.17,2.16) and (8.17,1.84) .. (m4-h103);
\draw[draw=black!45,line width=0.85pt] (m4-h103) .. controls (8.17,1.36) and (8.17,1.04) .. (m4-h101);
\draw[draw=black!45,line width=0.85pt] (m4-h101) .. controls (8.17,0.56) and (8.17,0.24) .. (m4-h100);
\draw[draw=black!45,line width=0.85pt] (m4-h102) .. controls (8.7,1.848) and (8.7,0.552) .. (m4-h100);
\node[dot] at (m4-h022) {};
\node[anchor=east,inner sep=0pt] at (-0.8,14.7) {$(u_{\xi},\xi)$};
\node[dot] at (m4-h023) {};
\node[anchor=east,inner sep=0pt] at (-0.8,13.9) {$(u_{\xi},\xi^2)$};
\node[dot] at (m4-h021) {};
\node[anchor=east,inner sep=0pt] at (-0.8,13.1) {$(u_{\xi},1)$};
\node[dot] at (m4-h020) {};
\node[anchor=east,inner sep=0pt] at (-0.8,12.3) {$(u_{\xi},0)$};
\node[dot] at (m4-h032) {};
\node[anchor=east,inner sep=0pt] at (-0.8,10.6) {$(u_{\xi^2},\xi)$};
\node[dot] at (m4-h033) {};
\node[anchor=east,inner sep=0pt] at (-0.8,9.8) {$(u_{\xi^2},\xi^2)$};
\node[dot] at (m4-h031) {};
\node[anchor=east,inner sep=0pt] at (-0.8,9) {$(u_{\xi^2},1)$};
\node[dot] at (m4-h030) {};
\node[anchor=east,inner sep=0pt] at (-0.8,8.2) {$(u_{\xi^2},0)$};
\node[dot] at (m4-h012) {};
\node[anchor=east,inner sep=0pt] at (-0.8,6.5) {$(u_{1},\xi)$};
\node[dot] at (m4-h013) {};
\node[anchor=east,inner sep=0pt] at (-0.8,5.7) {$(u_{1},\xi^2)$};
\node[dot] at (m4-h011) {};
\node[anchor=east,inner sep=0pt] at (-0.8,4.9) {$(u_{1},1)$};
\node[dot] at (m4-h010) {};
\node[anchor=east,inner sep=0pt] at (-0.8,4.1) {$(u_{1},0)$};
\node[dot] at (m4-h002) {};
\node[anchor=east,inner sep=0pt] at (-0.8,2.4) {$(u_{0},\xi)$};
\node[dot] at (m4-h003) {};
\node[anchor=east,inner sep=0pt] at (-0.8,1.6) {$(u_{0},\xi^2)$};
\node[dot] at (m4-h001) {};
\node[anchor=east,inner sep=0pt] at (-0.8,0.8) {$(u_{0},1)$};
\node[dot] at (m4-h000) {};
\node[anchor=east,inner sep=0pt] at (-0.8,0) {$(u_{0},0)$};
\node[dot] at (m4-h122) {};
\node[anchor=west,inner sep=0pt] at (8.8,14.7) {$(v_{\xi},\xi)$};
\node[dot] at (m4-h123) {};
\node[anchor=west,inner sep=0pt] at (8.8,13.9) {$(v_{\xi},\xi^2)$};
\node[dot] at (m4-h121) {};
\node[anchor=west,inner sep=0pt] at (8.8,13.1) {$(v_{\xi},1)$};
\node[dot] at (m4-h120) {};
\node[anchor=west,inner sep=0pt] at (8.8,12.3) {$(v_{\xi},0)$};
\node[dot] at (m4-h132) {};
\node[anchor=west,inner sep=0pt] at (8.8,10.6) {$(v_{\xi^2},\xi)$};
\node[dot] at (m4-h133) {};
\node[anchor=west,inner sep=0pt] at (8.8,9.8) {$(v_{\xi^2},\xi^2)$};
\node[dot] at (m4-h131) {};
\node[anchor=west,inner sep=0pt] at (8.8,9) {$(v_{\xi^2},1)$};
\node[dot] at (m4-h130) {};
\node[anchor=west,inner sep=0pt] at (8.8,8.2) {$(v_{\xi^2},0)$};
\node[dot] at (m4-h112) {};
\node[anchor=west,inner sep=0pt] at (8.8,6.5) {$(v_{1},\xi)$};
\node[dot] at (m4-h113) {};
\node[anchor=west,inner sep=0pt] at (8.8,5.7) {$(v_{1},\xi^2)$};
\node[dot] at (m4-h111) {};
\node[anchor=west,inner sep=0pt] at (8.8,4.9) {$(v_{1},1)$};
\node[dot] at (m4-h110) {};
\node[anchor=west,inner sep=0pt] at (8.8,4.1) {$(v_{1},0)$};
\node[dot] at (m4-h102) {};
\node[anchor=west,inner sep=0pt] at (8.8,2.4) {$(v_{0},\xi)$};
\node[dot] at (m4-h103) {};
\node[anchor=west,inner sep=0pt] at (8.8,1.6) {$(v_{0},\xi^2)$};
\node[dot] at (m4-h101) {};
\node[anchor=west,inner sep=0pt] at (8.8,0.8) {$(v_{0},1)$};
\node[dot] at (m4-h100) {};
\node[anchor=west,inner sep=0pt] at (8.8,0) {$(v_{0},0)$};
\end{scope}
\end{scope}
\end{tikzpicture}
\caption{The base graph $G_4$ and the four perfect matchings of cross edges in the MMS graph $H_4$.
In this case, we can construct $\mathbb{F}_4$ using the irreducible polynomial $x^2+x+1$ in $\mathbb{F}_2[x]$, that is, $\mathbb{F}_4=\mathbb{F}_2[x]/(x^2+x+1)$.
The elements of $\mathbb{F}_4$ can be represented by $0,1,x,x+1$.
Put $\xi=x+1$.
Then $\xi^2=x$, $\xi^3=1$, $\mathbb{F}_4=\{\xi,\xi^2,1,0\}$, $X=\{1,\xi^2\}=\{1,x\}$ and $X'=\{\xi,\xi^3\}=\{x+1,1\}$.
The union of cross edges in (d)--(g) is the set of cross edges of $H_4$.
The edges lifted from semiedges appear repeatedly in (d)--(g).
}
\label{fig:q=4}
\end{center}
\vspace{-0.5cm}
\end{figure}

\end{document}